\documentclass{svproc}
\usepackage{bm}
\usepackage{xcolor}
\usepackage{amsfonts} 
\usepackage{amsmath,amssymb} 
\usepackage{subcaption}
\usepackage{graphicx,float}
\usepackage{calc}
\usepackage[export]{adjustbox}
\usepackage{siunitx}

\newcommand{\TauFH}{\mathcal{T}_{F}}      
\newcommand{\TauAH}{\mathcal{T}_{A}}

\newcommand{\spanOp}[1]{\big\{#1\big\}}

\DeclareMathOperator*{\supp}{supp}
\DeclareMathOperator{\trace}{tr}

\usepackage{url}

\begin{document}
	\mainmatter              
	\title{AeroSPEED: a high order acoustic solver for aeroacoustic applications}
	\titlerunning{AeroSPEED: a high order acoustic solver for aeroacoustic applications}  
	%
 \author{}
	\author{Alberto Artoni\inst{1}, Paola F. Antonietti\inst{1}, Roberto Corradi\inst{2}, Ilario Mazzieri\inst{1}, Nicola Parolini\inst{1}, Daniele Rocchi\inst{2}, Paolo Schito\inst{2}, Francesco F. Semeraro\inst{2}}

	\authorrunning{Alberto Artoni et al.} 
	%
	%
	\institute{MOX, Dipartimento di Matematica, Politecnico di Milano, \\
		Piazza Leonardo da Vinci 32, Milano, 20133, Italy\\
		\email{alberto.artoni@polimi.it},
		\and
		Dipartimento di Meccanica, Politecnico di Milano, \\
		Via La Masa 1, Milano, 20156 , Italy}
	
	\maketitle              
	\begin{abstract}
We propose AeroSPEED, a solver based on the Spectral Element Method (SEM) that solves the aeroacoustic Lighthill's wave equation.  First, the fluid solution is computed employing a cell centered Finite Volume method. Then, AeroSPEED maps the sound source coming from the flow solution onto the acoustic grid, where finally the Lighthill's wave equation is solved. 
An ad-hoc projection strategy is adopted to apply the flow source term in the acoustic solver.
A model problem with a manufactured solution and the Noise Box test case are used as benchmark for the acoustic problem.  We studied the noise generated by the complex flow field around tandem cylinders as a relevant aeroacoustic application. AeroSPEED is an effective and accurate solver for both acoustics and aeroacoustic problems.
	\end{abstract}
	\section{Introduction}
%
Aeroacoustic is the field of acoustics that studies the noise induced by flows.
Due to the different scales involved in the flow and acoustics, usually aeroacoustics problems are posed in a segregated manner \cite{Kaltenbacher2020}.
First, a flow problem is solved and, in this work, the open-source finite volume library OpenFOAM \cite{OpenFOAM} is adopted.
Then, the flow solution is post-processed to generate the sound source term of the Lighthill's wave equation. With this purpose, we have developed AeroSPEED \cite{Artoni2023,SPEED}, a spectral element based solver that maps the computed sound source term onto the acoustic grid and solves the inhomogeneous Lighthill wave equation. We validate the open-source acoustic solver AeroSPEED on a model problem based on a manufactured solution, comparing it with COMSOL \cite{COMSOL}, a commercial software based on Lagrangian Finite Element Method (FEM). As an additional acoustic test case, we considered a geometry representing a simplified cockpit of a car (Noise Box). Next, we apply our aeroacoustic solver AeroSPEED to study the noise induced by the turbulent flow around two tandem circular cylinders. 
\begin{figure}
    \centering
\includegraphics[width=0.55\textwidth]{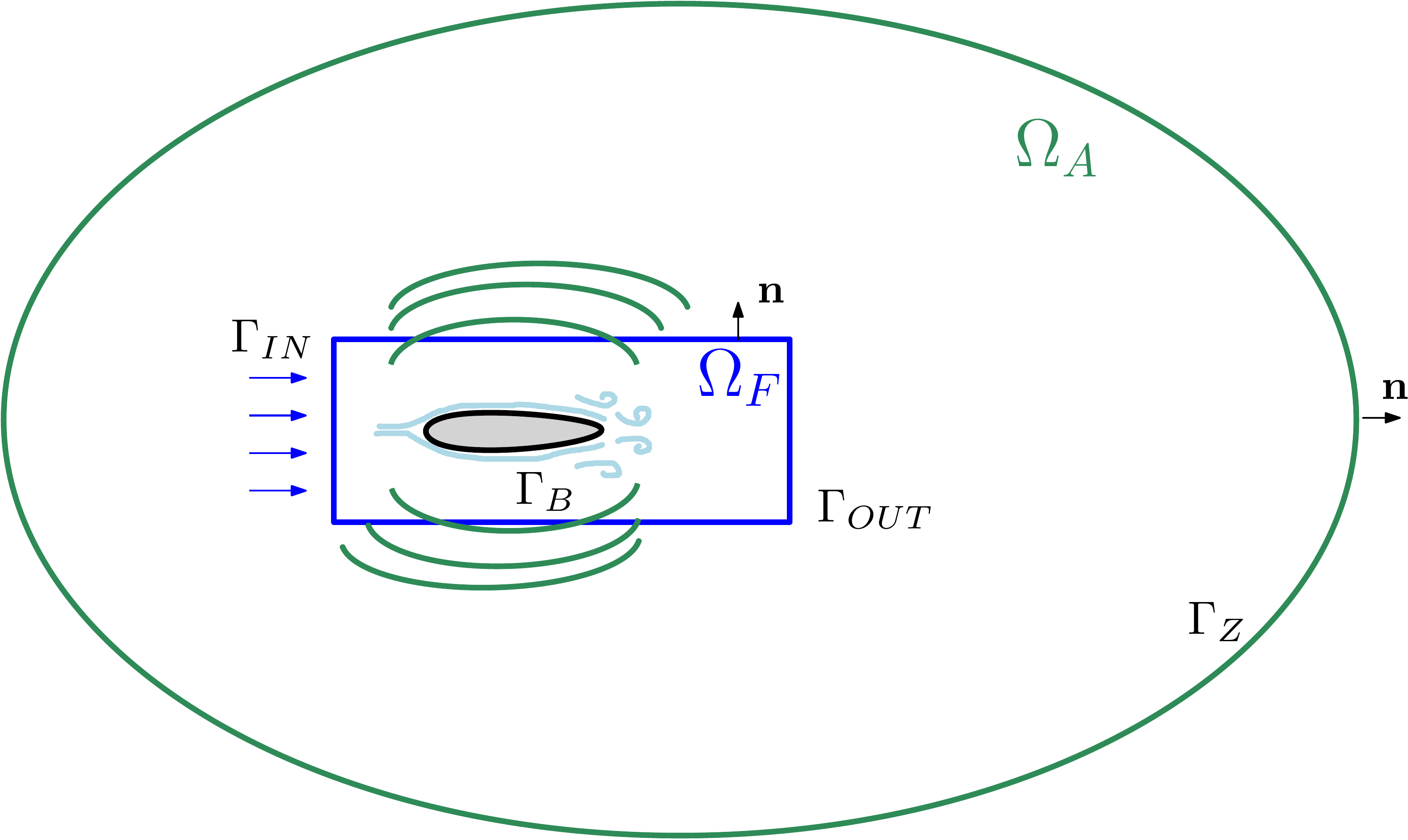}
    \caption{Domain of the aeroacoustic problem. $\Omega_F$ is the fluid domain, while $\Omega_A$ is the acoustic domain.}
    \label{fig:dominioAeroacustica}
\end{figure}
\section{The coupled aeroacoustic model problem}
	Given two open, bounded domains $\Omega_F \subseteq \Omega_A \subseteq \mathbb{R}^d$, having sufficiently regular boundaries $\partial \Omega_F$ and $\partial \Omega_A$ respectively (see Fig.~\ref{fig:dominioAeroacustica}), we consider the flow on a rigid body at high Reynolds and low Mach numbers. We are interested in the noise generated by an incompressible and acoustically compact flow, meaning that the feedback between the acoustic pressure and the hydrodynamic pressure can be neglected. Hence, the coupled aeroacoustic problem can be posed in a segregated manner.  The segregated approach considers the following sequence of problems.
	\subsubsection{Flow problem.} 
 The fluid flow is governed by the incompressible Navier Stokes equations: for $t \in (0,T]$, find the velocity field $\mathbf{u}(\mathbf{x},t):\Omega_F \times (0,T]\rightarrow \mathbb{R}^3$ and the pressure field $p(\mathbf{x},t):\Omega_F \times (0,T]\rightarrow \mathbb{R}$ such that
	\begin{equation}\label{eq:NavierStokes}
		\begin{aligned} 
			\frac{\partial\mathbf{u}}{\partial t} + \nabla\cdot(\mathbf{u}\otimes\mathbf{u}) - \nabla\cdot(\nu \nabla \mathbf{u}) + \nabla \left(\frac{p}{\rho_0}\right) &= 0, \quad \text{in } \Omega_F\times(0,T],\\
			\nabla\cdot\mathbf{u} &= 0, \quad \text{in } \Omega_F\times(0,T],\\
			\mathbf{u} &= \textbf{0}, \quad \text{on } \Gamma_B, \\
			\mathbf{u} &= \mathbf{g}, \quad \text{on } \Gamma_{IN}, \\
			\nu \nabla \mathbf{u}\cdot\mathbf{n} - p\mathbf{n} &= \textbf{0}, \quad \text{on } \Gamma_{OUT},
		\end{aligned}
	\end{equation}
	with initial condition $\mathbf{u}(\mathbf{x},0) = \textbf{0}$, and where $\mathbf{n}$ is the outward unit normal vector to $\partial\Omega_F$, $\nu$ is the kinematic viscosity, $\rho_0$ is the fluid density and $\mathbf{g}$ is the inlet Dirichlet datum. Moreover, we suppose the boundary of the fluid domain to be decomposed in the inlet boundary $\Gamma_{IN}$, the outlet boundary $\Gamma_{OUT}$ and the boundary $\Gamma_B$, such that $\partial\Omega_F = \Gamma_{IN}\cup\Gamma_{OUT}\cup \Gamma_{B}$.
	\subsubsection{Aeroacoustic source.} Based on the Lighthill analogy \cite{Kaltenbacher2020}, from the fluid velocity $\mathbf{u}$ we compute the Lighthill's tensor as $\mathbf{T} = \rho_0 \mathbf{u}\otimes\mathbf{u}$. The Lighthill's tensor has support only on the fluid domain $\Omega_F \subseteq \Omega_A$, and it depends only on the solution $\mathbf{u}$ of problem \eqref{eq:NavierStokes}. The Lighthill's tensor represents the sound source term and the coupling term between the flow problem \eqref{eq:NavierStokes} and the acoustic problem \eqref{eq:LighthillWaveContCont}.
	\subsubsection{Acoustic problem.} We consider in $\Omega_A$  the following non-homogeneous acoustic wave equation:
	for $t \in (0,T]$, find the density field $\rho(\mathbf{x},t):\Omega_A \times (0,T]\rightarrow \mathbb{R}$ such that
	\begin{equation}\label{eq:LighthillWaveContCont}
		\begin{aligned} 
			\frac{\partial^2\rho}{\partial t^2} -c_0^2 \Delta \rho = f, &\qquad \text{in } \Omega_A \times (0,T),   \\
			c_0^2 \frac{\partial\rho}{\partial\mathbf{n}} = 0, & \qquad \text{on } \Gamma_B\times (0,T),\\
			\frac{1}{\rho_0}\frac{\partial \rho}{\partial\mathbf{n}} = -\frac{1}{Z}\frac{\partial\rho}{\partial t}, & \qquad \text{on } \Gamma_Z\times (0,T),
		\end{aligned}
	\end{equation}	
	with initial conditions $\displaystyle \rho(\mathbf{x},0) = 0, \frac{\partial\rho}{\partial t}(\mathbf{x},0) = 0$, where $c_0$ is the sound speed and $Z$ is the impedance of an external wall. The boundary $\partial\Omega_A$ has been split as  $\partial\Omega_A = \Gamma_Z \cup \Gamma_B$ where $\Gamma_B$ is the body boundary where we set a sound hard boundary condition, while $\Gamma_Z$ is the external boundary where we apply an impedance boundary condition. We recall that if $Z$ is the characteristic impedance, i.e. $Z=\rho_0 c_0$, we have a non-reflective boundary condition, which is necessary for free-field wave propagation problems.
	When dealing with aeroacoustic problems the sound source is $f = \nabla\cdot \nabla\cdot \mathbf{T}$, obtaining the so called Lighthill's wave equation.
\section{Numerical scheme}
We introduce the spectral element method for the spatial discretization of \eqref{eq:LighthillWaveContCont} with a generic source term $f$, highlighting the aeroacoustic case where $f = \nabla\cdot \nabla\cdot \mathbf{T}$. 
\subsection{Spectral element dicretization} \label{sec:SEMapprox}
Given the acoustic domain $\Omega_A$, we introduce a conforming decomposition $\TauAH$ made by hexaedral elements $\kappa_A$ and we denote the characteristic mesh size as $\displaystyle h_A = \max_{\kappa_A \in \TauAH} h_{\kappa_A}$, being $h_{\kappa_A}$ the diameter of the element $\kappa_A$.
Let $\widehat{\kappa}$ be the reference element $\widehat{\kappa}=[-1,1]^3$, and assume that for any hexaedral element $\kappa_{A} \in \TauAH$ there exists a suitable trilinear invertible map $\bm{\theta}_{\kappa_{A}}:\widehat{\kappa}\rightarrow \kappa_{A}$, such that its Jacobian $\mathbf{J}_{\kappa_A}$ is positive.
We now introduce the following finite-dimensional space: 
\begin{equation}V_A = \left\{ v \in C^0(\overline{\Omega}_A)\cap H^1(\Omega_A): v|_{\kappa_{A}}\circ \bm{\theta}_{\kappa_A}^{-1}\in \mathbb{Q}_r(\widehat{\kappa}), \forall \kappa_{A} \in \TauAH\right\},
\end{equation}
where  $\mathbb{Q}_r(\widehat{\kappa})$ is the space of  polynomials of degree less than or equal to $r\geq 1$ in each coordinate direction.
Next, for any $u,w \in V_A$, we define the following bilinear form by means of the Gauss-Legendre-Lobatto (GLL) quadrature rule:
\begin{equation} \label{eq:GLLquad}
	(u,w) _{\kappa_{A}}^{NI}  = \sum_{i,j,k=0}^{r}  u(\bm{\theta}_{\kappa_{A}}(\bm{\xi}_{i,j,k}^{GLL}))w(\bm{\theta}_{\kappa_{A}}(\bm{\xi}_{i,j,k}^{GLL})) \omega_{i,j,k}^{GLL} |\det(\mathbf{J}_{\kappa_A})| \approx (u,w) _{\kappa_{A}},
\end{equation}
and we denote with
\begin{equation*}
	(u,w)_{\TauAH}^{NI} = \sum_{\kappa_A \in \TauAH} (u,w)_{\kappa_{A}}^{NI} \quad \forall \, u,w \in V_A,
\end{equation*}
where $\bm{\xi}^{GLL}$ are the GLL quadrature nodes, $\omega^{GLL}$ the corresponding weights defined on $\widehat{\kappa}$ and $NI$ stands for numerical integration. 
\subsection{Discretization of the acoustic problem}
\subsubsection{Weak formulation.}
We derive the weak formulation of the inhomogeneous wave equation in \eqref{eq:LighthillWaveContCont}:
\textit{for $t \in (0;T]$, find $\rho(\mathbf{x},t)\in H^1(\Omega_A)$ such that $\forall w \in H^1(\Omega_A)$:} 
\begin{align}
	\left(\frac{\partial^2\rho}{\partial t^2}, w\right)_{\Omega_A} + c_0^2(\nabla \rho ,\nabla w)_{\Omega_A}  + \frac{\rho_0 c_0^2}{Z} \int_{\Gamma_Z}\frac{\partial \rho}{\partial t} w\ ds  =  \mathcal{L}(w),
	\label{eq:BilinearContinuosProblem}  
\end{align}
with initial conditions $\rho = \displaystyle\frac{\partial \rho}{\partial t} = 0$ in $\Omega_A \times \{0\}$, and $\mathcal{L}(w)$ is a suitable linear operator.
When dealing with general inhomogeneous acoustic problems, the operator is \begin{equation}\mathcal{L}(w) = (f,w)_{\Omega_A}, \end{equation} while when dealing with Lighthill's wave equation the term $ (\nabla\cdot\nabla\cdot\mathbf{T}, w)_{\Omega_A}$ is usually integrated by parts, and we have that 
\begin{equation}\mathcal{L}(w) = -(\nabla\cdot\mathbf{T}, \nabla w)_{\Omega_A},\end{equation} see for instance \cite{Artoni2023}.
\subsubsection{Semi-discrete spectral element formulation.} For the sake of simplicity, we assume that $\partial\Omega_A = \Gamma_N$, and hence $\Gamma_Z = \emptyset$. The semi-discrete spectral element formulation of problem \eqref{eq:BilinearContinuosProblem} with numerical integration (SEM-NI) reads:
\textit{for any time $t\in(0;T]$ find $\rho_h \in V_A$ such that:}
\begin{equation}
	(\frac{\partial^2\rho_h}{\partial t^2},w_h)_{\TauAH}^{NI}+  c_0^2(\nabla \rho_h ,\nabla w_h)_{\TauAH}^{NI}  = \mathcal{L}_h(w_h) \quad \forall w_h \in V_A,
	\label{eq:LighthillWaveDiscreteCont}
\end{equation} 
with $\rho_h = \displaystyle\frac{\partial \rho_h}{\partial t} = 0$ in $\Omega_A \times \{ 0 \}$, where $\mathcal{L}_h(w) = (f,w)_{\TauAH}^{NI}$ for the purely acoustic case, while $\mathcal{L}_h(w) = -(\nabla\cdot\mathbf{T}, \nabla w)_{\TauAH}^{NI}$ for the aeroacoustic case.
\subsubsection{Computation of the aeroacoustic source term.}
Let $\TauFH$ be a polyhedral decomposition of the fluid domain $\Omega_F$ and let \begin{equation} V_F = \left\{ v_F \in L^2(\Omega_F):v_F|_{\kappa_F} \in \mathbb{P}^0(\kappa_F), \forall \kappa_F \in \TauFH \right\}
\end{equation} be the space of piecewise discontinuous functions. The sound source term $\nabla\cdot(\rho_0 \mathbf{u}\otimes\mathbf{u})$ is computed as a post-process of the numerical solution of problem \eqref{eq:NavierStokes} on the fluid grid $\TauFH$ via a Gauss discretization with a linear reconstruction. Namely, given the velocity $\mathbf{u}_h^k \in V_F$ at time $t^k$, we compute the cell value of the sound source term on the cell $\kappa_F \in \TauFH$ as:
\begin{multline} \label{eq:Lighthill}
	\frac{1}{|\kappa_F|} \int_{\kappa_F} \rho_0\nabla\cdot(\mathbf{u}_h^k\otimes\mathbf{u}_h^k) d\mathbf{x} = \frac{1}{|\kappa_F|} \int_{\partial\kappa_F} \rho_0 \mathbf{u}_h^k (\mathbf{u}_h^k\cdot\mathbf{n}_{\mathcal{F}}) d s\approx \\ \frac{1}{|\kappa_F|} \sum_{\mathcal{F} \in \partial\kappa_F} \rho_0\mathbf{u}_\mathcal{F} (\mathbf{u}_\mathcal{F}\cdot\mathbf{n}_{\mathcal{F}} )|\mathcal{F}|,
\end{multline}
where $\mathbf{u}_{\mathcal{F}} = \mathbf{u}(\mathbf{x}_{\mathcal{F}})$, $\mathbf{n}_\mathcal{F}$ is the unit normal to the face  $\mathcal{F}\in \partial\kappa_F$ and where we applied a mid-point quadrature rule using the mid-point $\mathbf{x}_{\mathcal{F}}$ of face $\mathcal{F}$. The value of the velocity at the face centre $\mathbf{x}_{\mathcal{F}}$ is computed with a linear interpolation.

\subsubsection{Projection of the aeroacoustic source term.}
Let $q_F \in V_F$ be a function defined on the fluid grid $\TauFH$ such that $q_F = \sum_{i=1}^{N_F} \widehat{q}_{F,i} \phi_{F,i}$, where  $\spanOp{\phi_{F,i}}_i^{N_F}$ is the set of $N_F$ basis functions associated to $V_F$, and $\widehat{q}_{F,i}$ are the corresponding expansion coefficients. 
We define the $L^2$-projection of the field $q_F \in V_F$ into $V_A$ as:\
\textit{find $q_A \in V_A$ s.t. } 
\begin{align}
	( q_{A}, \phi_{A,i})_{\TauAH} &= ( q_F,  \phi_{A,i})_{\TauAH}  \quad \forall \phi_{A,i} \in V_{A}, \label{def:ProjectionA}
\end{align}
where $q_A \in V_A$ is a function defined on the discrete acoustic space such that $q_A = \sum_{i=1}^{N_A} \widehat{q}_{A,i} \phi_{A,i}$, where  $\spanOp{\phi_{A,i}}_i^{N_A}$ is the set of $N_A$ basis functions, and $\widehat{q}_{A,i}$ are the corresponding expansion coefficients.
Since $q_F$ is a piecewise constant over $\TauFH$, namely $q_F \in V_F$, we recast problem  (\ref{def:ProjectionA}) as follows:
\begin{multline}  \label{eq:ProjectionProblemDiscrete}
	\sum_{\kappa_{A} \in \TauAH} ( q_{A}, \phi_{A,i})_{\kappa_A} = \sum_{\kappa_{A} \in \TauAH} \left( \sum_{\ell=1}^{N_F} \widehat{q}_{F,\ell} \phi_{F,\ell} , \phi_{A,i} \right)_{\kappa_A} = \\ \sum_{\kappa_{A} \in \TauAH} \sum_{\ell=1}^{N_F}  \widehat{q}_{F,\ell} ( 1, \phi_{A,i} )_{\kappa_A\cap \kappa_{F,\ell} }, 
\end{multline}
where we have used that $\kappa_{F,\ell} =  \supp(\phi_{F,\ell})$.
The discrete algebraic counterpart of  \eqref{eq:ProjectionProblemDiscrete} becomes
\begin{equation} \label{eq:AlgebraicProblemProjection}
	\mathbf{M}^{AA} \widehat{\mathbf{q}}_A = \mathbf{M}^{AF} \widehat{\mathbf{q}}_{F},
\end{equation}
where $\mathbf{M}^{AA} \in \mathbb{R}^{N_A\times N_A}$ with
\begin{equation}
    \mathbf{M}^{AA}_{i,j} = \sum_{\kappa_{A} \in   \TauAH} (\phi_{A,j}, \phi_{A,i})_{\kappa_A}, \qquad i,j=1,\dots,N_A,
\end{equation} 
is the full mass matrix and
$\mathbf{M}^{AF} \in \mathbb{R}^{N_A\times N_F}$  is defined as
\begin{equation}
	\mathbf{M}^{AF}_{i,\ell} =  \sum_{\kappa_{A} \in   \TauAH} \int_{\kappa_A \cap \kappa_{F,\ell}} \phi_{A,i} \ \text{d}\mathbf{x}, \qquad  i=1,\dots,N_A, \; \ell=1,\dots,N_F. \label{eq:rhsProj}
\end{equation}
The vector $\widehat{\mathbf{q}}_A$ in \eqref{eq:AlgebraicProblemProjection} collects the expansion coefficients of the projected acoustic field $q_{A}$, while $\widehat{\mathbf{q}}_F$ collects the expansion coefficients of the donor fluid field $q_{F}$.
Further details on the implementation of the projection method are given in \cite{Artoni2023}.
\subsubsection{Algebraic formulation of the semi-discrete problem.}
We introduce the mass and stiffness matrix $\mathbf{M},\mathbf{K} \in \mathbb{R}^{N_A\times N_A}$:
\begin{equation} \label{eq:massMatrix}
	\mathbf{M}_{i,j} = \sum_{\kappa_{A} \in \TauAH} ( \phi_{A,j},\phi_{A,i})^{NI}_{\kappa_A}, \qquad
	\mathbf{K}_{i,j} = \sum_{\kappa_{A} \in \TauAH} ( \nabla\phi_{A,j},\nabla\phi_{A,i})^{NI}_{\kappa_A},
\end{equation} 
for $i,j=1,\dots,N_A$. We remark that the mass matrix $\mathbf{M}$ computed with the quadrature formula in eq.~\eqref{eq:GLLquad} becomes diagonal and hence ${\mathbf{M} \neq \mathbf{M}^{AA}}$.
Since $\rho_h \in V_A$, we write $\displaystyle \rho_h = \sum_{i=1}^{N_A} \widehat{\rho}_{A,i} \phi_{A,i}$, where  $\spanOp{\phi_{A,i}}_i^{N_A}$ is the set of $N_A$ basis functions associated to $V_A$. Collecting the expansion coefficients $\widehat{\rho}_{A,i}$ into the vector $\bm{\rho}_h$, we obtain the following algebraic semi-discrete formulation:
\begin{equation} \label{eq:semiDiscreteAcoustics}
	\mathbf{M} \ddot{\bm{\rho}}_h + c_0^2\mathbf{K}\bm{\rho}_h = \mathbf{f},
\end{equation}
supplemented with the initial conditions $\bm{\rho}_h = \mathbf{0}$ and $\dot{\bm{\rho}}_h = \mathbf{0}$. \\
For an inhomogeneous acoustic problem we define 
\begin{equation}
	\mathbf{f}_{i} = \mathcal{L}_h(\phi_{A,i}) = \sum_{\kappa_{A} \in \TauAH} ( f,\phi_{A,i})^{NI}_{\kappa_A}, \quad i=1,\dots,N_A,
\end{equation}
while for the aeroacoustic problem we have that:
\begin{equation}
    \mathbf{f}_{i} = \Big[\sum_{\ell=x,y,z}\mathbf{C}^\ell\widehat{\mathbf{q}}_{A,\ell}\Big]_i, \quad i=1,\dots,N_A.
\end{equation}
For the aeroacoustic case in fact, given $\widehat{\mathbf{q}}_{A,\ell}$ solution of the projection problem \eqref{eq:ProjectionProblemDiscrete} with $\widehat{\mathbf{q}}_{F,\ell} = [\nabla\cdot\mathbf{T}]_\ell$, computed as described in Eq.~\eqref{eq:Lighthill}, with $l=x,y,z$ representing each component, we have that:
\begin{equation}
\begin{aligned}
    &\mathcal{L}_h(\phi_{A,i}) =
	\sum_{\ell=x,y,z} \left(\sum_{\kappa_A \in \TauAH} -(q_{A,\ell}, [\nabla \phi_{A,i}]_\ell)_{\kappa_A}^{NI} \right) =\\ &\sum_{\ell=x,y,z} \left( \sum_{\kappa_A \in \TauAH} - \Big( \sum_{j=1}^{N_A}  \widehat{q}_{A,\ell}\phi_{A,j},[\nabla \phi_{A,i}]_\ell \Big)^{NI}_{\kappa_A} \right) = \Big[\sum_{\ell=x,y,z}\mathbf{C}^\ell\widehat{\mathbf{q}}_{A,\ell}\Big]_i,
 \end{aligned}
\end{equation}
with $i=1,\dots,N_A,$ and where $\mathbf{C}^\ell$ is defined as
\begin{equation}
	\mathbf{C}^\ell_{i,j} = \sum_{\kappa_A \in \TauAH} (\phi_{A,j} , [\nabla \phi_{A,i}]_\ell)_{\kappa_A}^{NI}, \qquad \text{for } i,j=1,\dots,N_A.
\end{equation}
\subsubsection{Time discretization.} 
For the time discretization of problem \eqref{eq:semiDiscreteAcoustics} we divide the temporal interval $(0, T]$ into $N$ subintervals, such that $T = N\Delta t$, and we set $t^k = k\Delta t$, with $k=0,\dots,N-1$ and introduce the auxiliary variables $\displaystyle \bm{v}_h^k = \dot{\bm{\rho}}_h(t^k)$, $\bm{a}_h^k = \ddot{\bm{\rho}_h}(t^k)$.
Furthermore, since the mass matrix $\mathbf{M}$ is not singular, we can represent Eq.~\eqref{eq:semiDiscreteAcoustics} as:
\begin{equation}
	\ddot{\bm{\rho}_h} = \mathcal{A}(\bm{\rho}_h,t), \label{eq:fullyDiscreteacousticsPb}
\end{equation}
where $\mathcal{A}(\bm{\rho}_h,t) = \mathbf{M}^{-1}(\mathbf{f}-c_0^2\mathbf{K}\bm{\rho}_h)$.
We now employ the Newmark method to discretize  Eq.~\eqref{eq:fullyDiscreteacousticsPb}:
\begin{equation}
\begin{aligned} \label{eq:fullyDiscreteacousticsNw}
	\bm{\rho}^{k+1}_h &= \bm{\rho}^{k}_h + \Delta t \bm{v}^{k}_h +\Delta t^2 \Big(\beta_N \mathcal{A}^{k+1} + (\frac{1}{2}-\beta_N) \mathcal{A}^{k}\Big) ,\\
	\bm{v}^{k+1}_h &= \bm{v}^{k}_h + \Delta t \Big(\gamma_N \mathcal{A}^{k+1} + (1-\gamma_N) \mathcal{A}^{k} \Big), \\
\end{aligned}
\end{equation}
where $\displaystyle 0 \le \beta_N \le \frac{1}{2}$ and $0\le \gamma_N \le 1$ are parameters of the Newmark scheme, and where $\mathcal{A}^k = \mathcal{A}(\bm{\rho}^k_h,t^k)$.
%
%
%
\section{Curle Analogy}
\label{sec:Curle} By following the derivation in \cite{Curle1955}, we recall the Curle aeroacoustic solution for problem \eqref{eq:LighthillWaveContCont}, that will be considered for comparison in the numerical results presented in Section~\ref{sec:NumericalResults}. Given an observer located at $\mathbf{x}$ at the time $t$, a volume $V$ and a body $B \subset V$, we have that:
\begin{multline}
p(\mathbf{x},t) = \int_{-\infty}^{+\infty} \int_{V} \frac{\partial^2 T_{ij}(\mathbf{y},\tau)}{\partial x_i \partial x_j} G(\mathbf{x},t|\mathbf{y},\tau) 	\textrm{d}\mathbf{y}\textrm{d}\tau \\ + \int_{-\infty}^{+\infty}\int_{\partial B} \left( p(\mathbf{y},\tau)\frac{\partial G(\mathbf{x},t|\mathbf{y},\tau)}{\partial \mathbf{n}} -  G(\mathbf{x},t|\mathbf{y},\tau)\frac{\partial p(\mathbf{y},\tau)}{\partial \mathbf{n}}\right)d\partial B \textrm{d}\tau, \label{eq:Kirchoff}
\end{multline}
where $G$ is a suitable Green function, $V$ is the control volume, $\mathbf{n}$ is the outward unit normal to the boundary $\partial B$.
We denote with $\mathbf{r} = \mathbf{x} - \mathbf{y}$, being $r$ its modulus. We choose as Green function $\displaystyle G(\mathbf{x},t|\mathbf{y},\tau) = \frac{1}{4\pi r} \delta(\tau - t + \frac{r}{c_0})$ in eq.~\eqref{eq:Kirchoff} to obtain:
	\begin{equation} \label{eq:KirchoffIntegralLighthill}
		\begin{aligned}
			p(\mathbf{x},t)=& \frac{1}{4\pi}  \int_V \frac{1}{r}\Big[\frac{\partial^2 T_{ij}}{\partial x_i \partial x_j}\Big]  d\mathbf{y}  \\+
			& \int_{\partial B} \frac{1}{4 \pi r}  \left[ \left(\frac{1}{c_0}\frac{\partial p}{\partial t} +
			\frac{p}{r}\right)\frac{\mathbf{r}}{r}\cdot\mathbf{n}
			- \frac{\partial p}{\partial \mathbf{n}} \right] d\partial B, \end{aligned}
	\end{equation}
 where $[\cdot]$ means that the function has to be evaluated at the retarded time $\displaystyle t-\frac{r}{c_0}$.
Next, we perform the following simplifications (see \cite{LibAcoustics} for details). First,  the volume term containing the Lighthill tensor $\mathbf{T}$ is neglected.
Then the retarded times are neglected. This assumption is reasonable if the considered sound sources are compact, that means if the characteristic length of the emitting object $D$ is smaller then the characteristic length $\lambda$ of the acoustic wave, namely if $D \ll \lambda$.
	Furthermore, since the object is considered acoustically rigid, i.e. $\displaystyle \left.\frac{\partial p}{\partial \mathbf{n}}\right\rvert_{\partial B} = 0$, we have that:
	\begin{align}
		p(\mathbf{x},t) & = \int_{\partial B} \frac{1}{4 \pi r}  \left(\left(\frac{1}{c_0}\frac{\partial p}{\partial t} + \frac{p}{r}\right)\frac{\mathbf{r}}{r}\cdot \mathbf{n} \right) d\partial B. \label{eq:CurleOF}
	\end{align}
	By neglecting the viscous forces and considering $\displaystyle \mathbf{F} = \int_{\partial B} p \mathbf{n}$, we have that:
	\begin{align}
		p(\mathbf{x},t) & = \frac{1}{4 \pi}\frac{\mathbf{r}}{r^2}\cdot\left(\frac{\mathbf{F}}{r} + \frac{1}{c_0}\frac{\partial\mathbf{F}}{\partial t}\right).
	\end{align}
\section{Numerical results for acoustic problems}
We consider the inhomogeneous wave equation described in \eqref{eq:LighthillWaveContCont} and we compare our software AeroSPEED based on the spectral element approximation introduced in Section~\ref{sec:SEMapprox} with the commercial software COMSOL \cite{COMSOL} based on the Lagrangian FEM.
\subsection{Verification test case}
As a first test case, we consider a simple model problem where we verify the performance of both AeroSPEED and COMSOL, in terms of accuracy and computational efficiency.
\subsubsection{Acoustic Setup.}
Given the manufactured solution \begin{equation}
    u_{ex}=\sin(\pi t)\sin(4\pi(x-1)(y-1)(z-1))\sin(4\pi xyz),
\end{equation}
we solve an inhomogeneous wave equation on the cube $\Omega_A = (0,1)^3$, with Neumann boundary conditions on $\Gamma_B = \partial \Omega_A$, with $c_0=\SI{1}{\meter\per\second}$. We employ a very fine time step $\Delta t=1\times10^{-6}\SI{}{\second}$ and we use a second order Newmark scheme with parameters $\gamma_N=0.5$ and $\beta_N=0.25$. We solve the test case both in AeroSPEED and in COMSOL, changing the refinement of the acoustic grid and the polynomial degree of the underlying polynomial approximation and we compute the error $E_2=||u_{ex}-u_{h}||_2$ at the final time $T=0.5\SI{}{\second}$.
\subsubsection{Acoustic Results.}
\begin{figure}
	\centering
 \begin{subfigure}{0.48\textwidth}
  \includegraphics[width=\textwidth]    {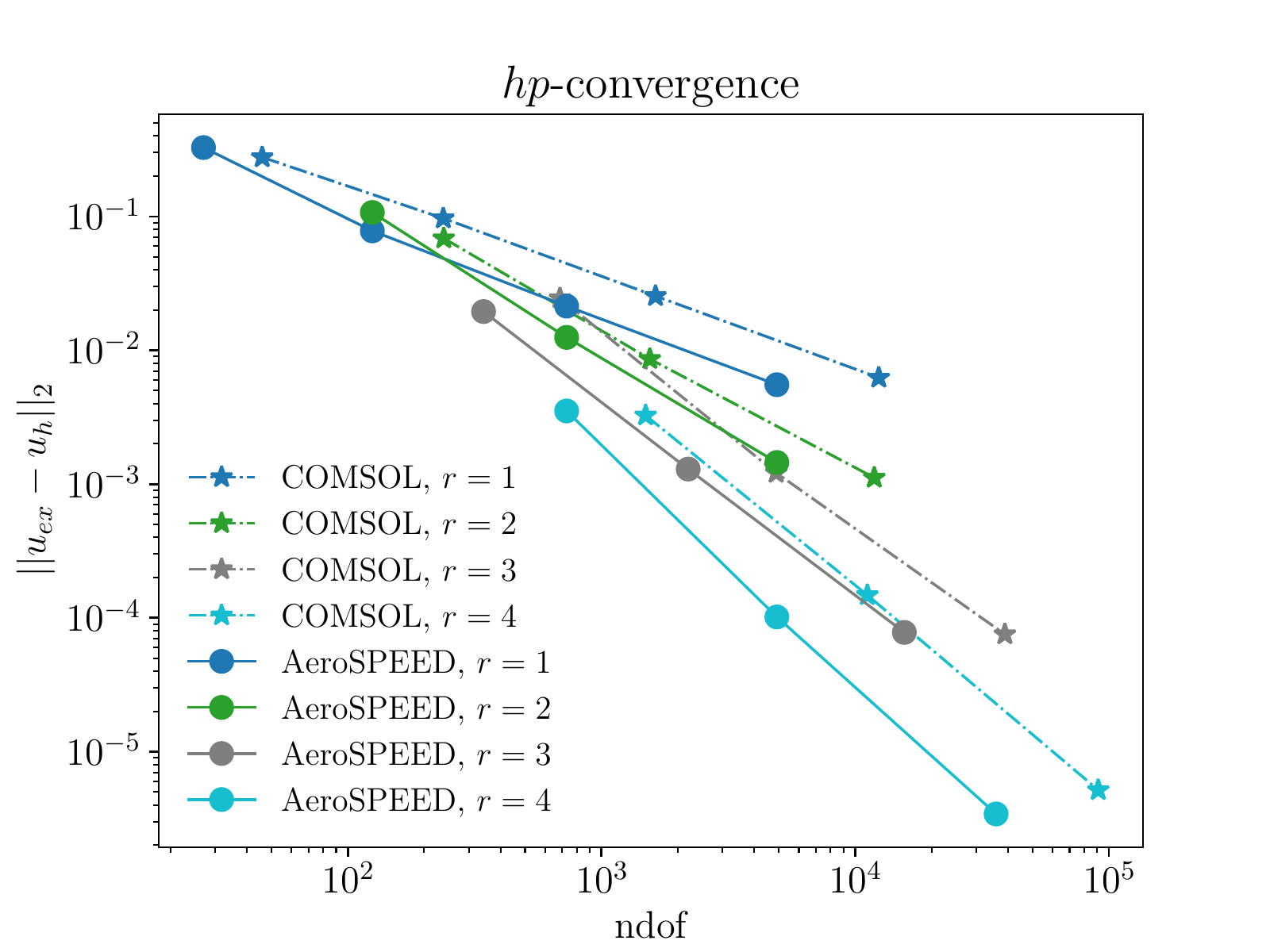} 
 \caption{}
 \end{subfigure}
  \begin{subfigure}{0.48\textwidth}
	\includegraphics[width=\textwidth]   {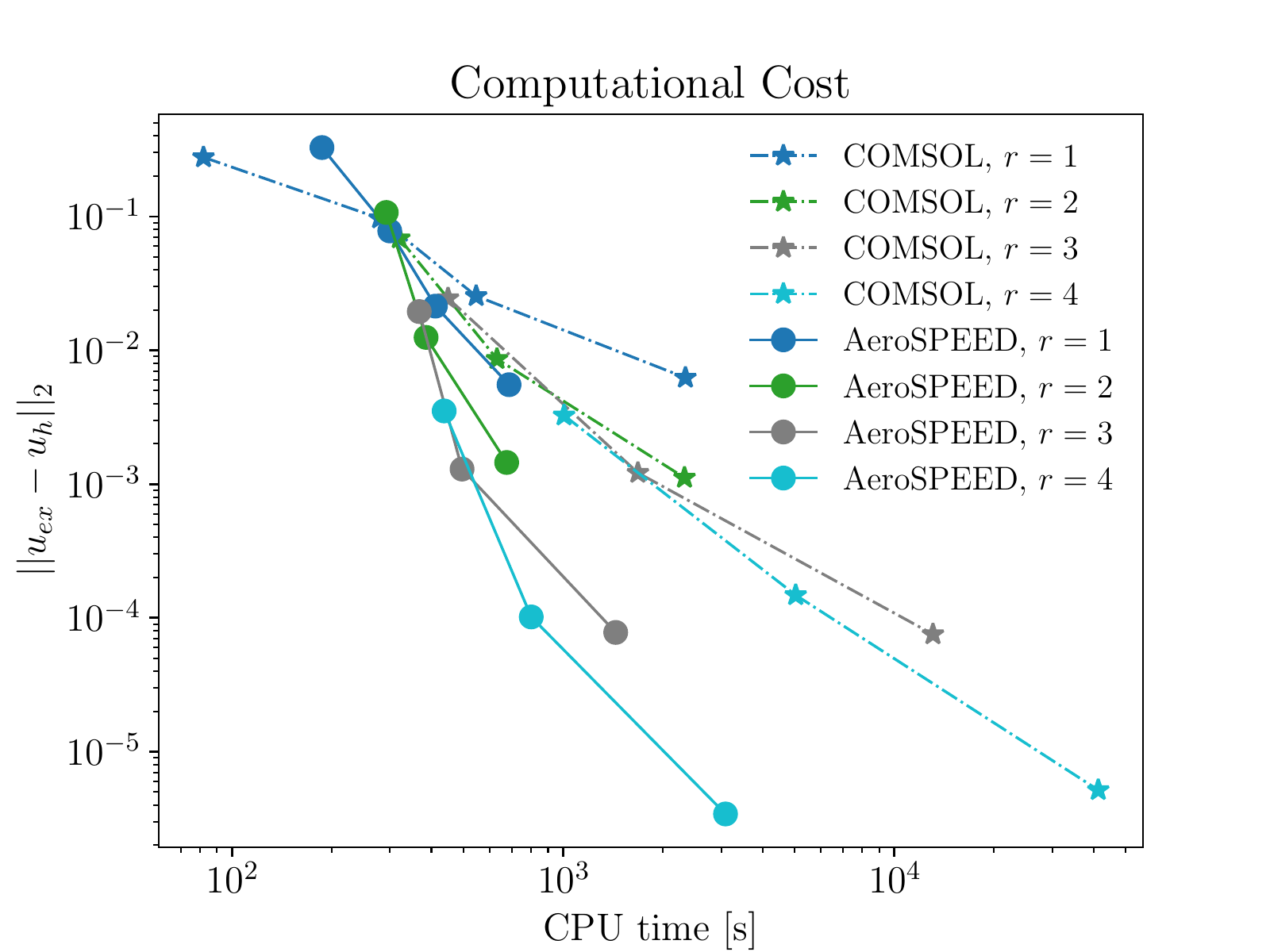}
	\caption{}
 \end{subfigure}
 \caption{Comparison between the SEM solver AeroSPEED and the Lagrangian FEM solver COMSOL. (a) Computed errors vs number of degrees of freedom (ndof). (b) Computed error versus CPU time. The computational tests were performed on 4 cores Intel(R) Xeon(R) Gold 6226 CPU at 2.70GHz.}
	\label{fig:COMSOLvsSPEED}
\end{figure}

 We report in Fig.~\ref{fig:COMSOLvsSPEED} the computed errors versus the number of degrees of freedom (left) and the CPU time (right) obtained with the AeroSPEED and COMSOL solvers varying the polynomial approximation degree $r=1,2,3,4$ of a sequence of meshes with comparable granularity.
 The expected convergence rates are obtained for both the underlying SEM and FEM approximations, respectively. For a comparable number of degrees of freedom, the SEM approximation is more accurate and less expensive.
 The results reported in Fig.~\ref{fig:COMSOLvsSPEED} (right) clearly indicate that AeroSPEED is able to achieve the same error in a much shorter computational time.
 Moreover, as the underlying polynomial approximation degree increases, AeroSPEED becomes more and more efficient compared to COMSOL.
\subsection{Noise Box test case}
We consider a second test case to asses the capabilities of AeroSPEED in solving acoustic problems in a confined geometry and we compare the obtained numerical solution with the one provided by COMSOL.
\subsubsection{Acoustic Setup.}
\begin{figure}
\centering
\begin{subfigure}{0.5\textwidth}
    \includegraphics[width=\textwidth]{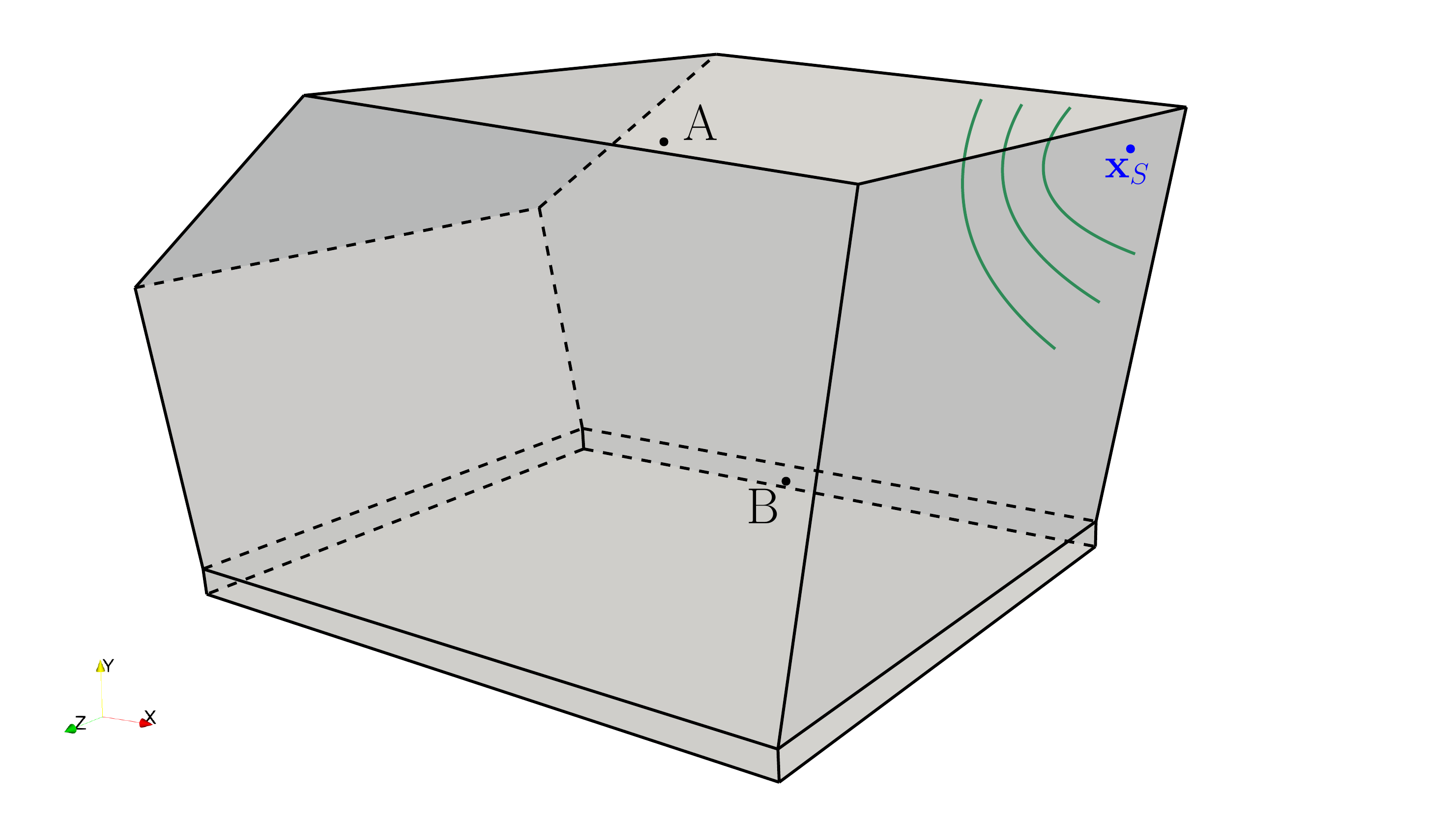}
    \caption{}
\end{subfigure}
\begin{subfigure}{0.25\textwidth}
    \includegraphics[width=\textwidth]{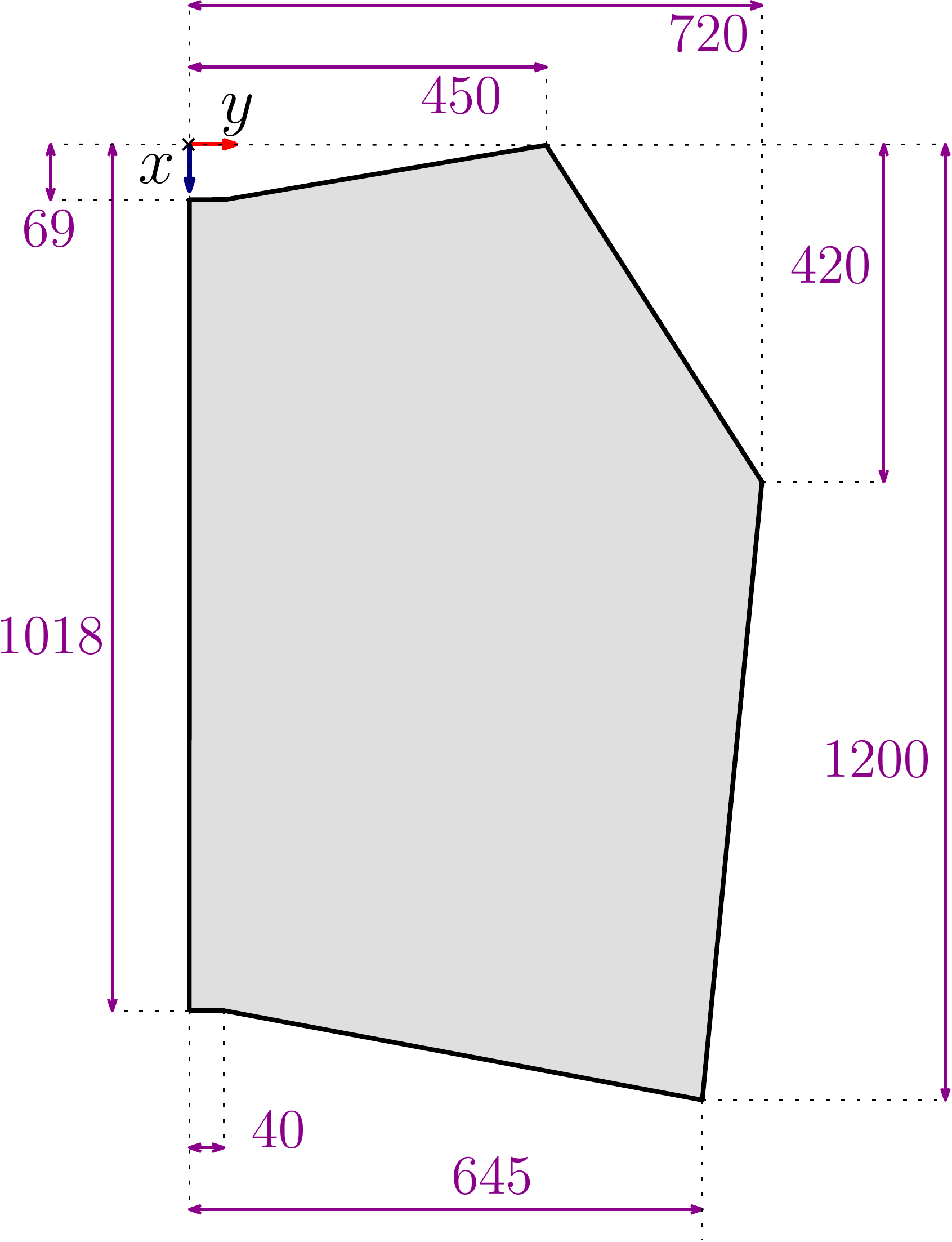}
    \caption{}
\end{subfigure}
\caption{
 (a) Three-dimensional view of the domain of the Noise Box. A and B are the positions of the selected microphones, where A = $(0.424,0.595,0.151)\ \SI{}{\meter}$ and B = $(0.9, 0.224, 0.528)\ \SI{}{\meter}$. (b) Quoted computational domain of the Noise Box. The spanwise length is $\SI{0.825}{\meter}$.  In the figure, units are expressed in millimeters.}
	\label{fig:NoiseBox}
\end{figure}

We consider a geometry that represents a simplified car cockpit, the so called Noise Box, see Fig.~\ref{fig:NoiseBox}, introduced in \cite{LingLiu2021}.
Each wall is modeled as a real wall (with both partially reflective and partially absorbing behaviour), by setting a wall impedance
of $Z= \SI{32206}{\pascal\second\per\meter}$.
As forcing term we consider a monopole sound source $f(\mathbf{x},t) = \delta(\mathbf{x} - \mathbf{x}_S)\sin(2 \pi f_0  t)$, where $\delta(\mathbf{x}-\mathbf{x}_S)$ is the Dirac delta centered in $\mathbf{x}_S = (1.15,0.595,0.065) \SI{}{\meter}$ and $f_0 = \SI{162}{\hertz}$. We set the density of air to be $\rho_0 = \SI{1.204}{\kilogram\per\meter\cubed}$ and the speed of sound $c_0 = \SI{343}{\meter\per\second} $. For the space discretization we set the polynomial degree $r=2$ and we fix $\Delta x =  \SI{0.04}{\meter}$. For the time discretization we set $\gamma_N = 0.5$ and $\beta_N = 0.25$, with $\Delta t = 5\times10^{-6} \SI{}{\second}$, with a final time of $T= \SI{0.5}{\second}$. We solve for the same setup both with COMSOL and AeroSPEED and we compare the two results.

\subsubsection{Acoustic Results.}
\begin{figure} 
	\centering
  \includegraphics[width=0.48\textwidth]    {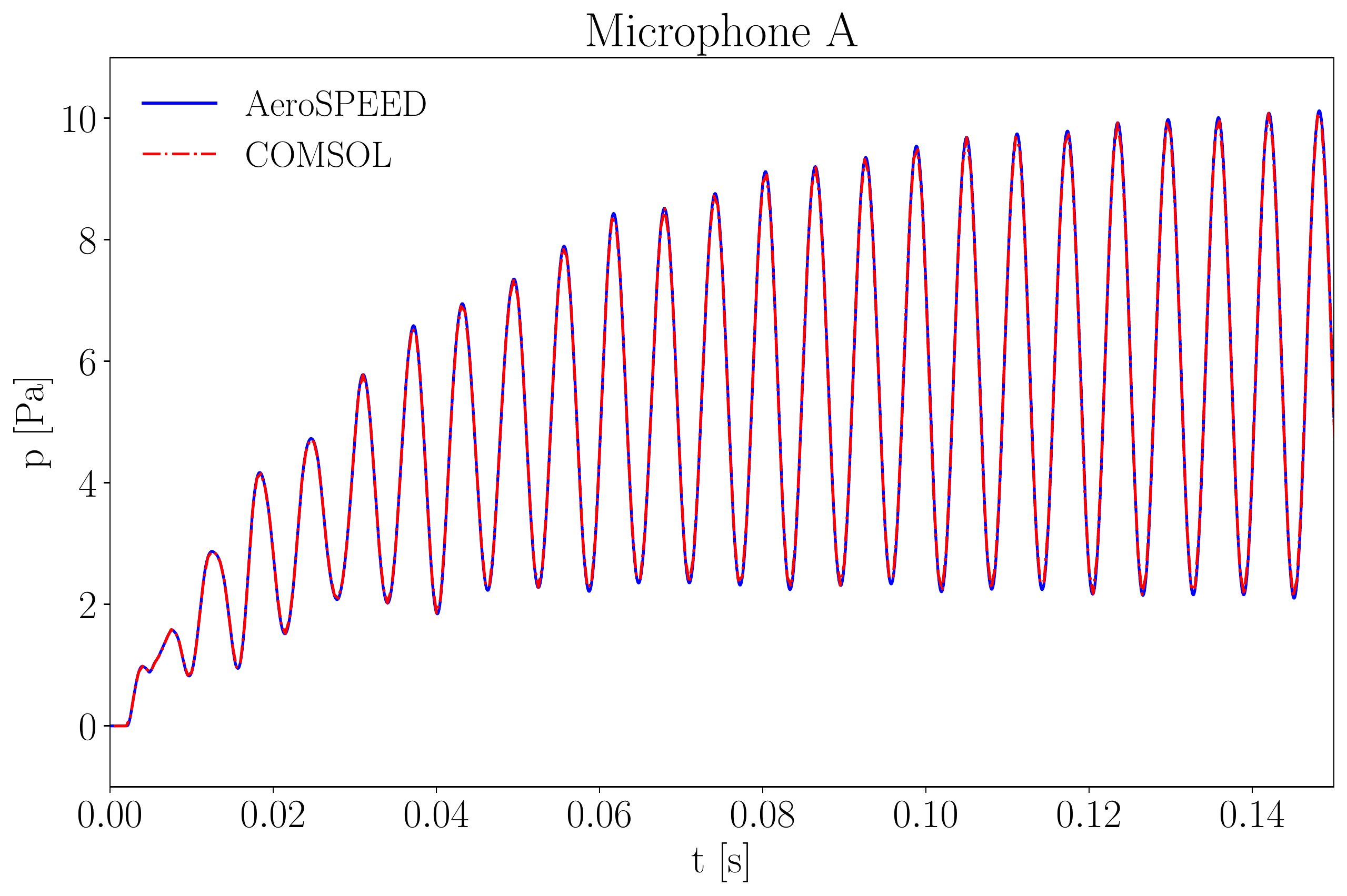} 
 \hfill
    \includegraphics[width=0.48\textwidth]    {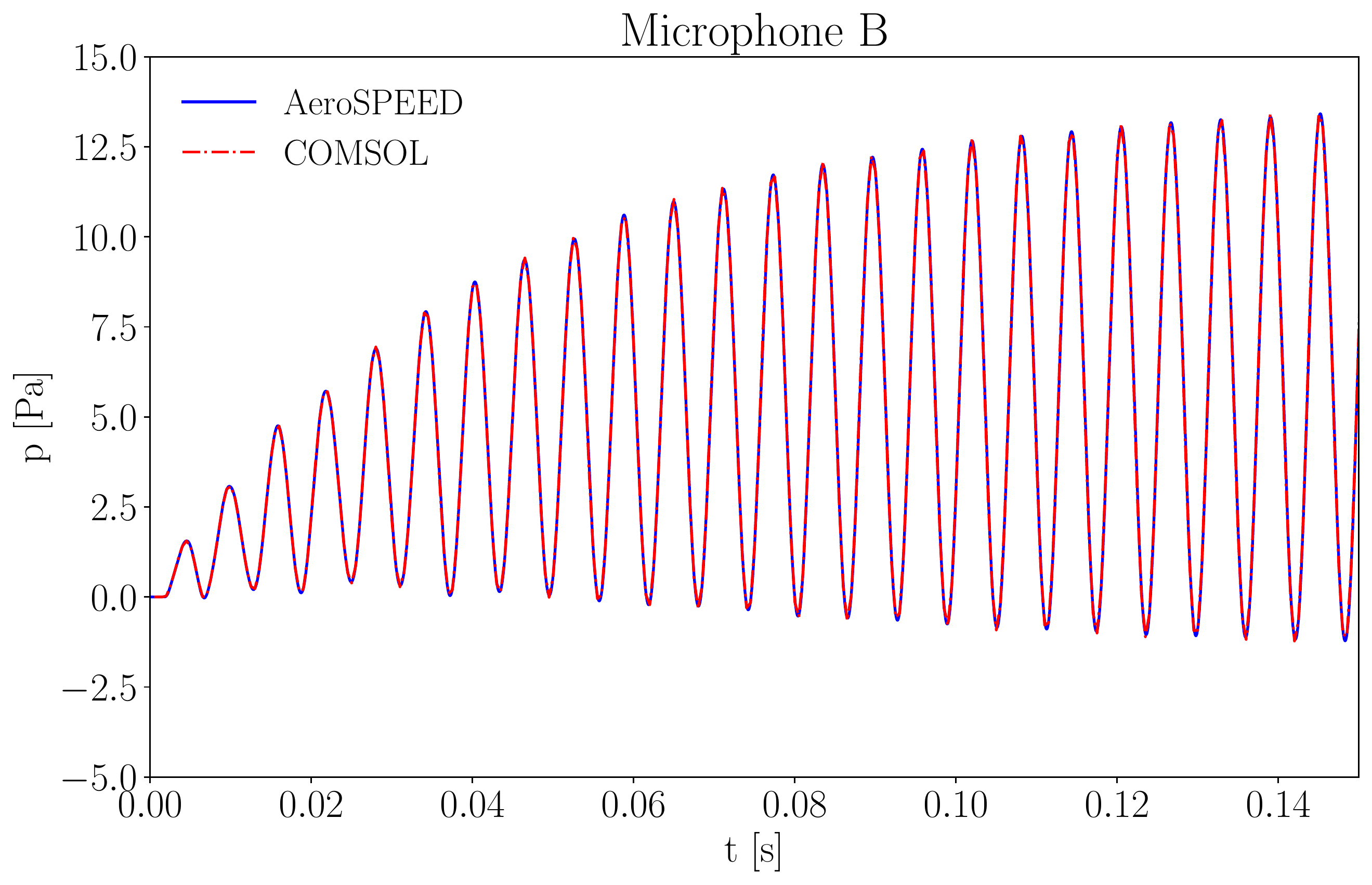} 
	\caption{Computed acoustic pressure measured by microphone A and B with both AeroSPEED and COMSOL.} 
	\label{fig:ProbesNoiseBox}
\end{figure}

\begin{figure}
\adjustbox{valign=t}{
\begin{subfigure}{0.28\textwidth}
    \includegraphics[width=\textwidth]{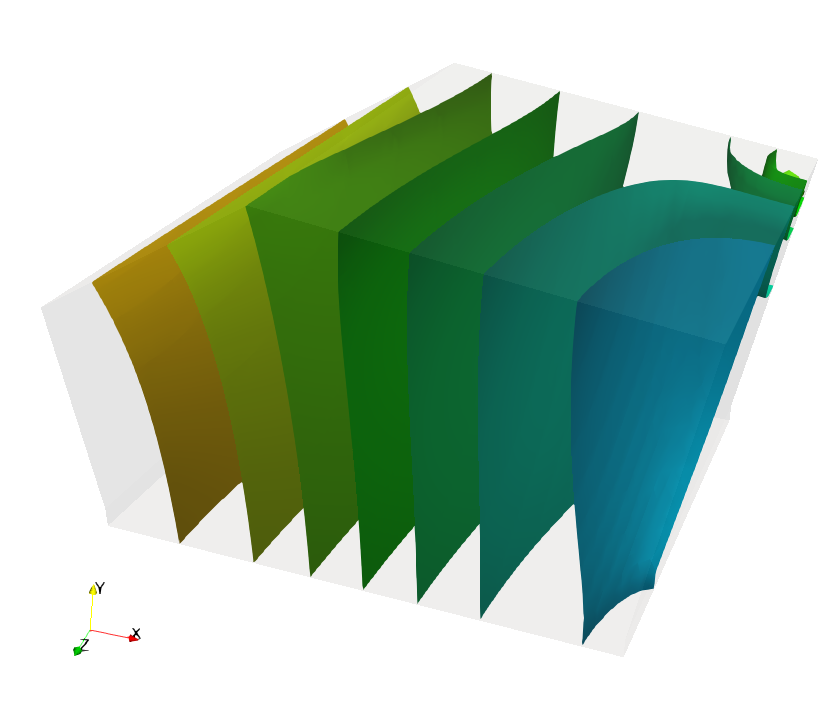}
    \caption{$t = \SI{0.495}{\second}$}
\end{subfigure} }
\adjustbox{valign=t}{
\begin{subfigure}{0.28\textwidth}
    \includegraphics[width=\textwidth]{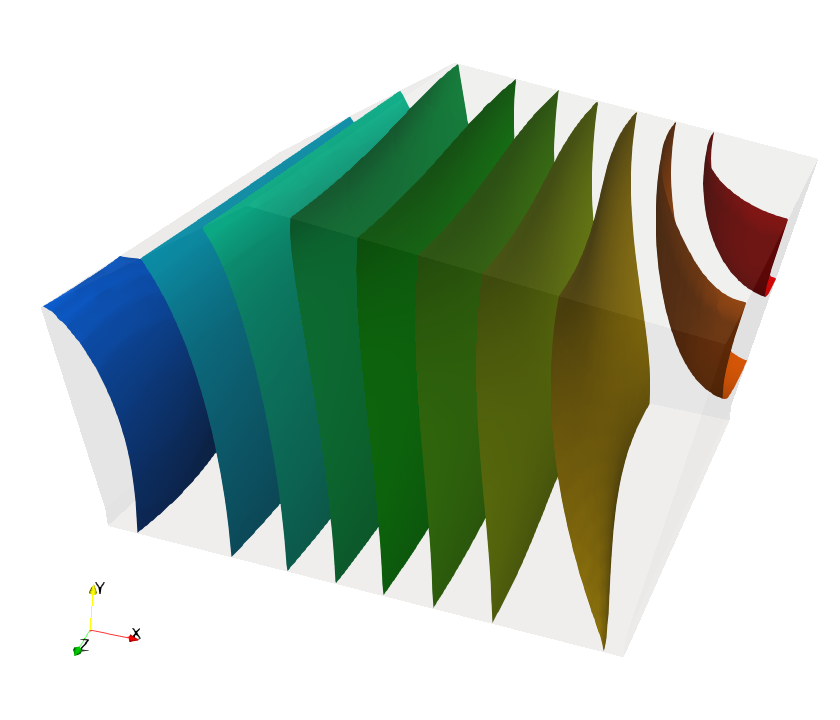}
        \caption{$t = \SI{0.496}{\second}$}
\end{subfigure} }
\adjustbox{valign=t}{
\begin{subfigure}{0.28\textwidth}
    \includegraphics[width=\textwidth]{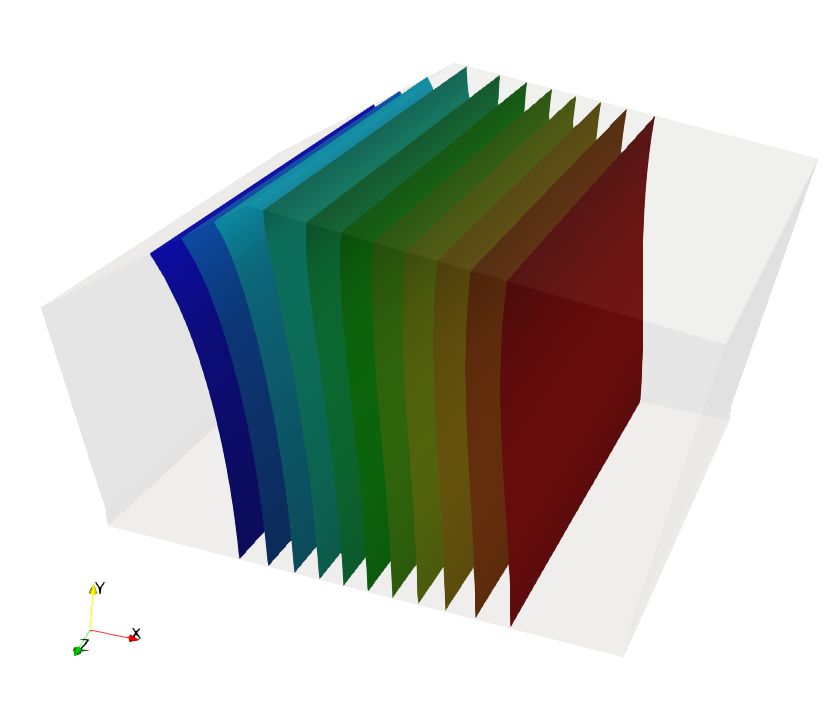}
        \caption{$t = \SI{0.497}{\second}$}
\end{subfigure}} \hfill
\adjustbox{valign=t}{
\begin{subfigure}{0.07\textwidth}
    \includegraphics[width=\textwidth]{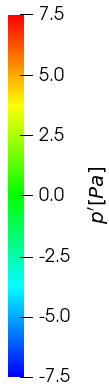}
\end{subfigure}}
    \caption{ Snapshots of the computed pressure fluctuations $p' = p - \overline{p}$, where $\overline{p}$ is the average pressure, inside the Noise Box, for $t=0.495, 0.496,  \SI{0.497}{\second}$. The selected contour levels are from $\SI{-7.5}{\pascal}$ to $\SI{7.5}{\pascal}$ with a step of $\SI{1.5}{\pascal}$.}
    \label{fig:PressureWaves}
\end{figure}

From the results reported in  Fig.~\ref{fig:ProbesNoiseBox}a we note the initial transient state, up to around $t \approx  \SI{0.05}{\second}$. The acoustic monopole is injecting energy in the system, that is not fully dissipated, up until $t \approx  \SI{0.1}{\second}$.
At that time, the system has reached a stationary regime, where the amount of energy dissipated by the system is balanced by the amount of energy injected. The numerical solution obtained with AeroSPEED perfectly matches the numerical solution obtained with COMSOL. 
In Fig.~\ref{fig:PressureWaves} we see the stationary pressure waves inside the Noise Box from different snapshots of the solution.
\section{Numerical results for an aeroacoustic application} \label{sec:NumericalResults}

We consider an aeroacoustic application, namely the noise generated by two tandem cylinders, a test case that have been subject of a dedicated workshop \cite{Lockard2011}.
The flow simulation has been performed with OpenFOAM \cite{OpenFOAM}, while the aeroacoustic coupling has been implemented in AeroSPEED \cite{Artoni2023}.
\subsection{Turbulent flow around a tandem cylinder}
Simulating the turbulent flow around  two tandem cylinders at high Reynolds number is a challenging problem due to the unsteadiness and complex flow structures to be captured.
The separation point on the front cylinder moves on the surface, generating a shear layer that rolls up forming a periodic vortex shedding that impinges on the rear cylinder. As result, a tonal and broadband noise are generated. Proper turbulence models are crucial to simulate at a reasonable computational cost such a complex physics. %
\begin{figure}
    \centering
    \adjustbox{valign=m}{
    \begin{subfigure}{0.47\textwidth}\includegraphics[width=\textwidth]{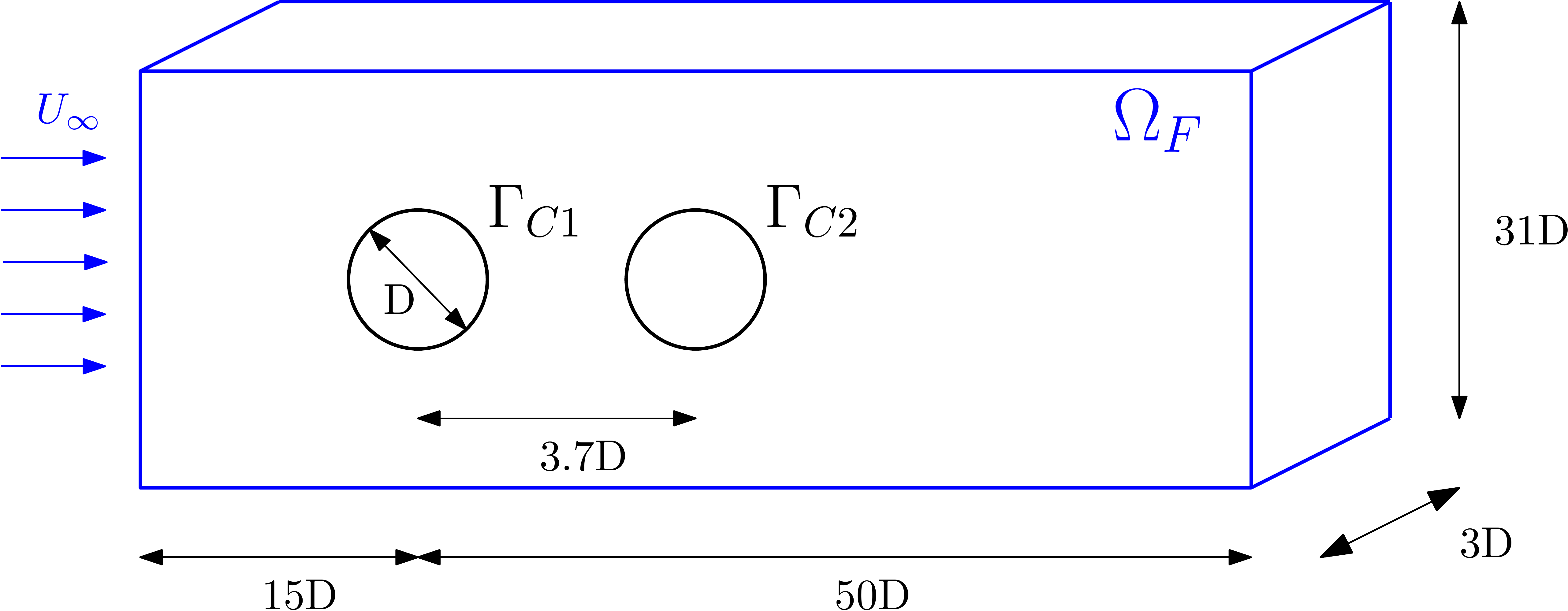}
    \caption{}
    \end{subfigure}}
    \adjustbox{valign=m}{
    \begin{subfigure}{0.47\textwidth}
    \includegraphics[width=\textwidth]{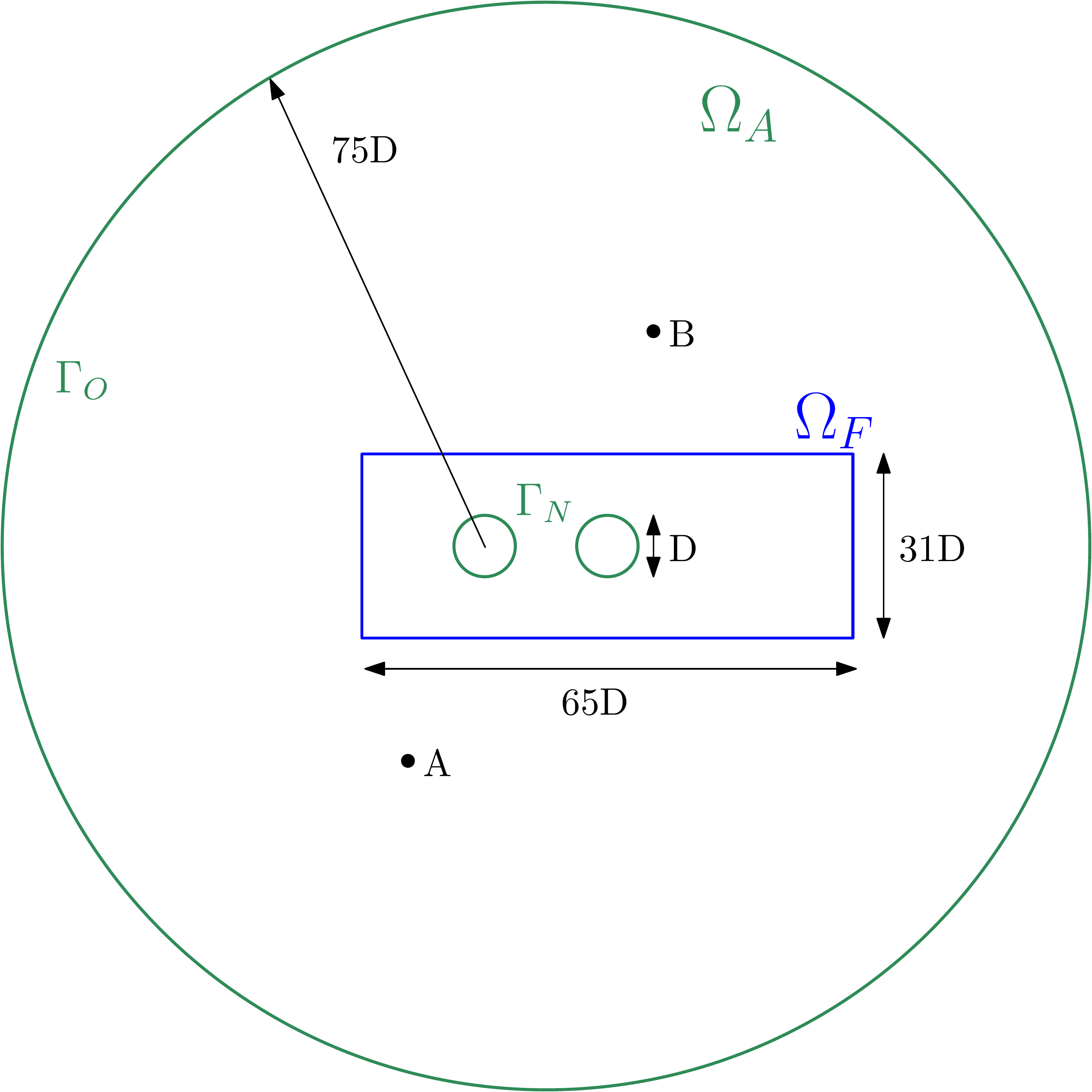}
    \caption{}
    \end{subfigure}}
    \caption{Tandem cylinders computational domain. (a) Fluid computational domain. (b) Aeroacoustics computational domain. The center of the computational domain is set at the center of the front cylinder. The points A, B are microphone probes located at A = (-8.33D,27.82D), and B = (9.11D, 32.49D), with $D= \SI{0.057}{m}$. The front cylinder is denoted by $\Gamma_{C1}$ and the rear cylinder with $\Gamma_{C2}$.}\label{fig:TandemComputationalDomains}
\end{figure}

\subsubsection{Fluid Setup.}
The two tandem cylinders problem configuration involves two cylinders of equal diameter $D = \SI{0.057}{\meter}$ aligned along the streamwise direction at a distance of $3.7D$. 
A sketch of the computational domain is reported in Fig. \ref{fig:TandemComputationalDomains}a.
At the inlet a fixed velocity of $U_\infty = \SI{44}{\meter\per\second}$ is set, corresponding to a Reynolds number $Re = 1.66\times10^5$. 
On $\Gamma_{C1}$ and $\Gamma_{C2}$, no slip conditions are imposed. At the outlet, a zero gradient condition is set. On the remaining boundaries, a symmetry condition is imposed. 
We choose a fixed time step of $\Delta t = 1.25\times10^{-5}\ \SI{}{\second}$, we set the final time to $T = \SI{0.35}{\second}$, and we employ a second order backward difference formula.
The height of the first cell near the wall corresponds to $y^+ \approx 30$, and proper wall functions are prescribed, see \cite{Spalding1961}.
Following \cite{Greschner2012}, we employ a DDES $k-\omega$SST model to simulate the turbulent flow, see for instance \cite{kOmegaDDES2011} for more details.
\subsubsection{Acoustic Setup.}
Since the acoustic problem can be considered bi-dimensional, we take as sound source only the average along the spanwise direction. A sketch of the computational domain for the aeroacoustic case is depicted in Fig. \ref{fig:TandemComputationalDomains}b. On the cylinders $\Gamma_{C1} \cup \Gamma_{C2} = \Gamma_N$ rigid wall boundary conditions are imposed, while at the far field $\Gamma_O$ absorbing boundary conditions are considered. The fluid sound source is mapped each four time steps, namely the computational time step for the Lighthill's wave equation is $\Delta t = 5\times10^{-5}\ \SI{}{\second}$. The chosen polynomial degree is $r=2$ and the spacing at the far field is $\Delta x \approx \SI{0.04}{\meter} $.
\subsubsection{Flow Validation.}
\begin{figure}
\adjustbox{valign=m}{	\begin{subfigure}{0.43\textwidth}
		\includegraphics[width=\textwidth]{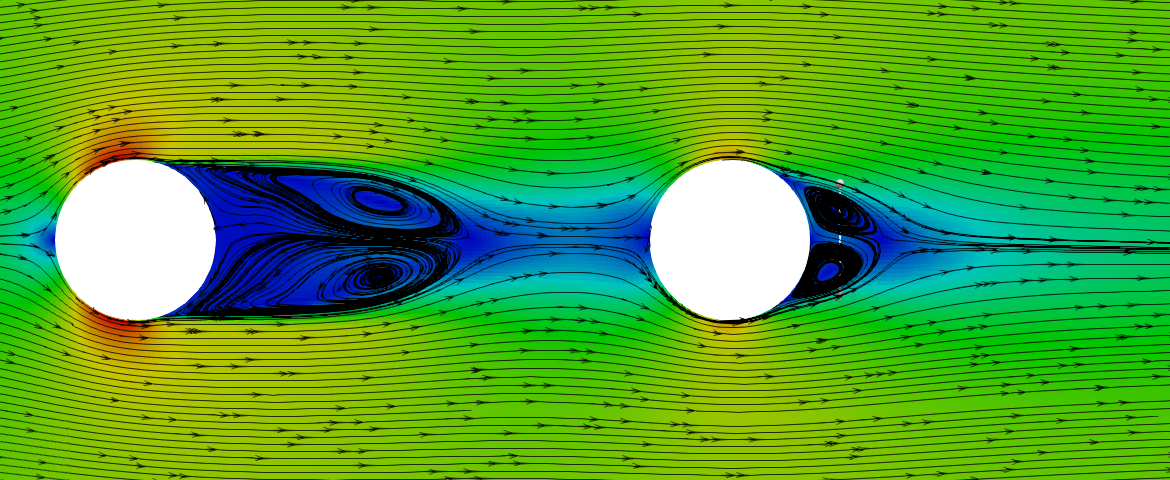} \caption{}
	\end{subfigure}
 }
 \adjustbox{valign=m}{
	\begin{subfigure}{0.43\textwidth}
		\includegraphics[width=\textwidth]{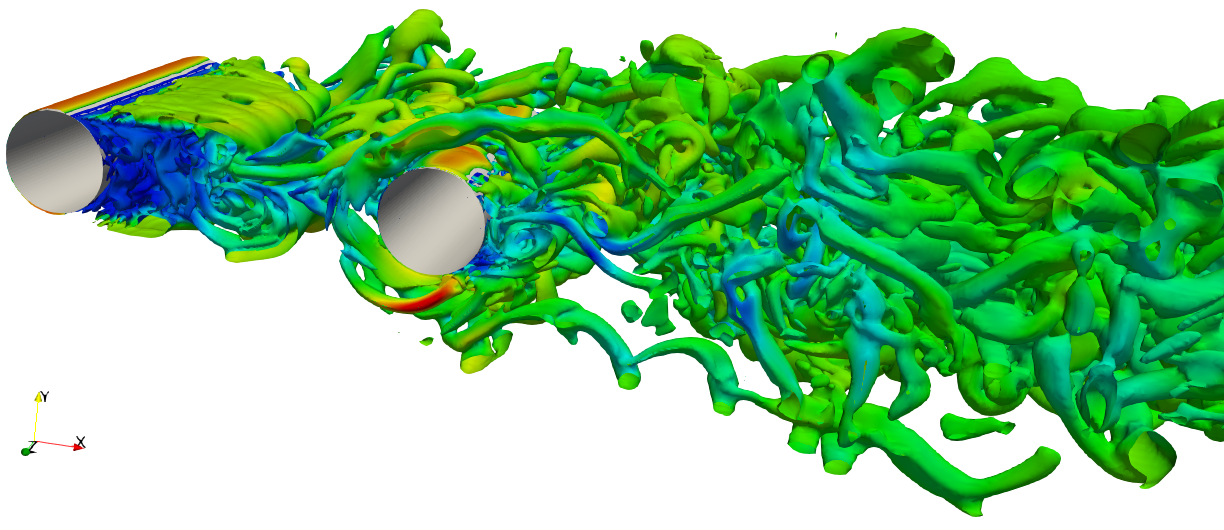}
 \caption{}
 \end{subfigure}	
 } \hfill
 \adjustbox{valign=m}{
 \begin{subfigure}{0.07\textwidth}
 \includegraphics[width=\textwidth]{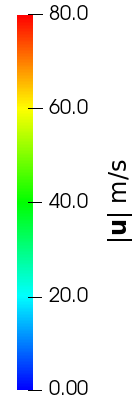}
    \end{subfigure}}
	\caption{
 (a) Average velocity magnitude $|\mathbf{u}|$ on the symmetry plane with streamlines. (b) Isosurfaces of Q=1000 at $t=\SI{0.35}{\second}$ colored with velocity magnitude.
 \label{fig:flowfield}}
\end{figure}

\begin{figure} 
\centering
    \includegraphics[valign=m,width=0.12\textwidth]{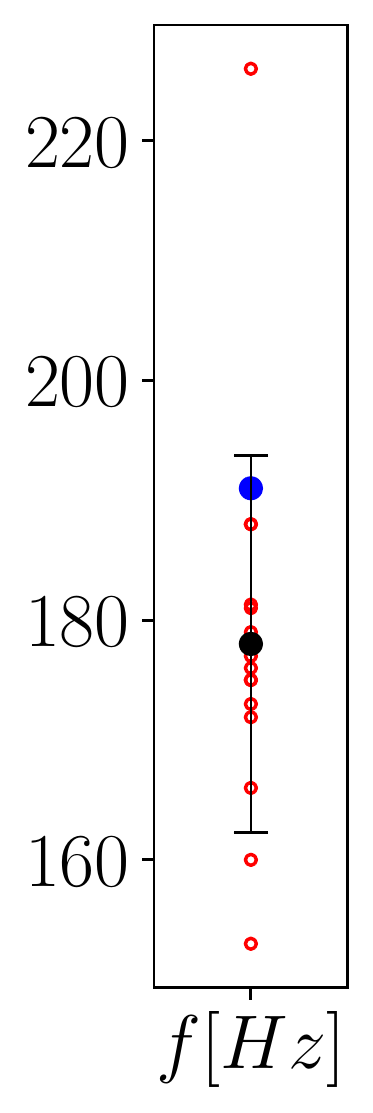}
    \includegraphics[valign=m,width=0.4\textwidth]{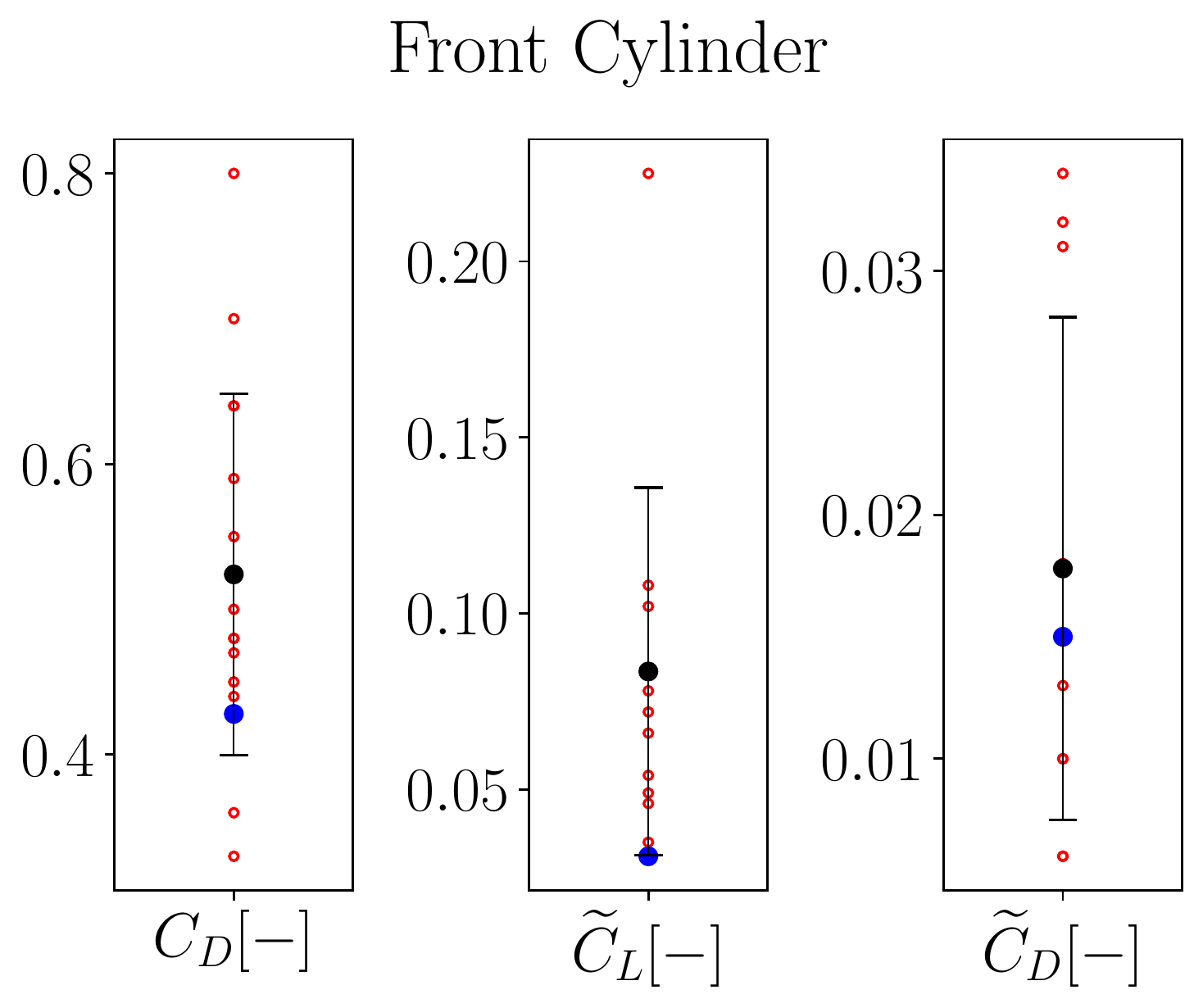}
    \includegraphics[valign=m,width=0.4\textwidth]{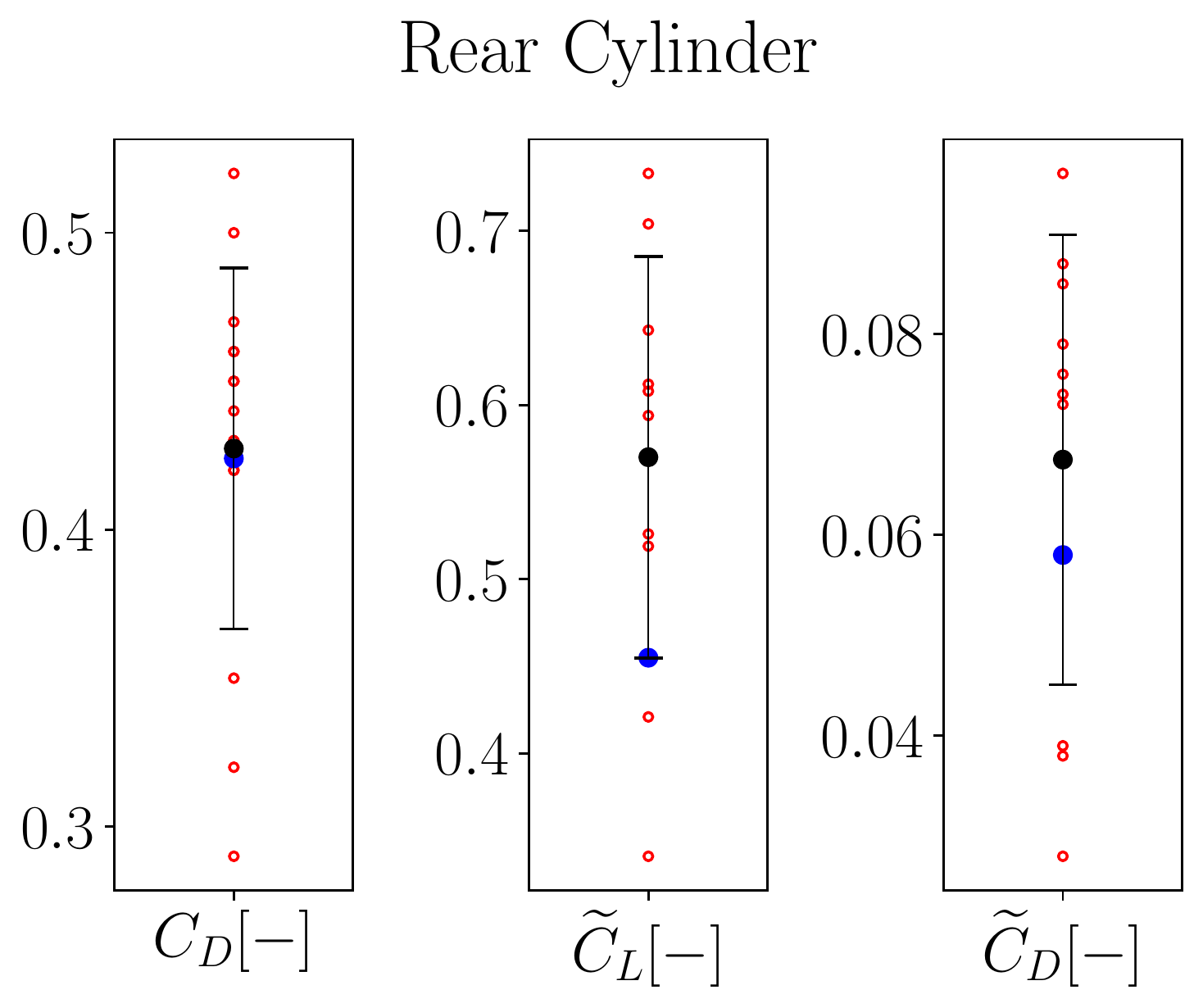}\\
    \includegraphics[width=0.55\textwidth]{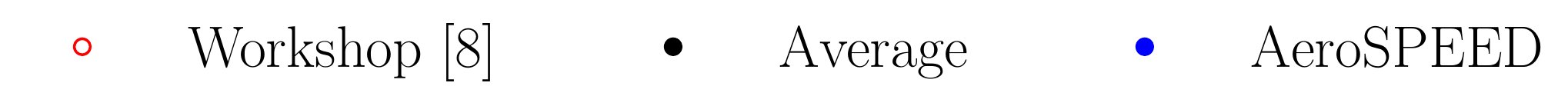} 
\caption{Comparison of the average drag forces and the $rms$ for lift and drag of the two cylinders. Data have been taken from \cite{Lockard2011}. The lift frequency is common between the two cylinders. The bar is centred on the mean value of the available literature data, and the length of the bar is one standard deviation. \label{fig:IntegralValues}}
\end{figure}

\begin{figure}
 \adjustbox{valign=m}{
 	\begin{subfigure}{0.11\textwidth}
  \includegraphics[width=\textwidth]{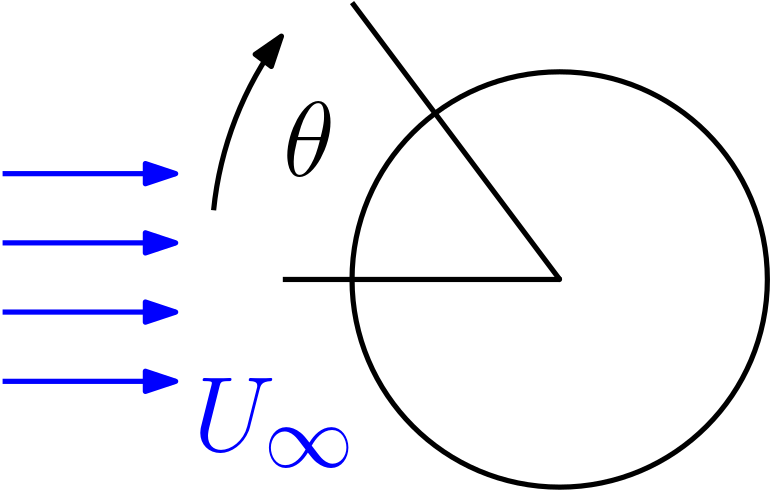}
    \end{subfigure}}
     \adjustbox{valign=m}{
	\begin{subfigure}{0.43\textwidth}
		\includegraphics[width=\textwidth]{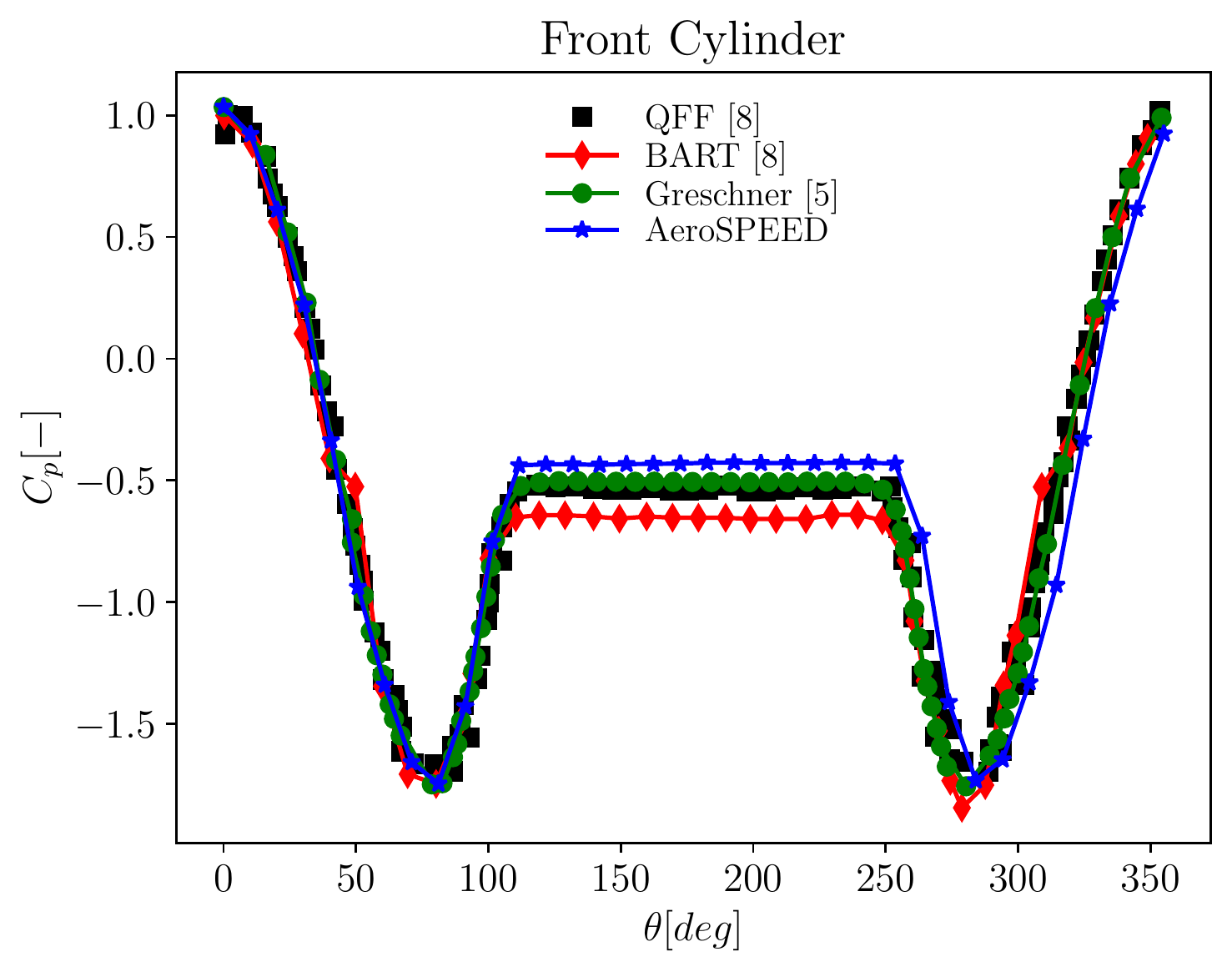}
	\end{subfigure}}
  \adjustbox{valign=m}{
	\begin{subfigure}{0.43\textwidth}
		\includegraphics[width=\textwidth]{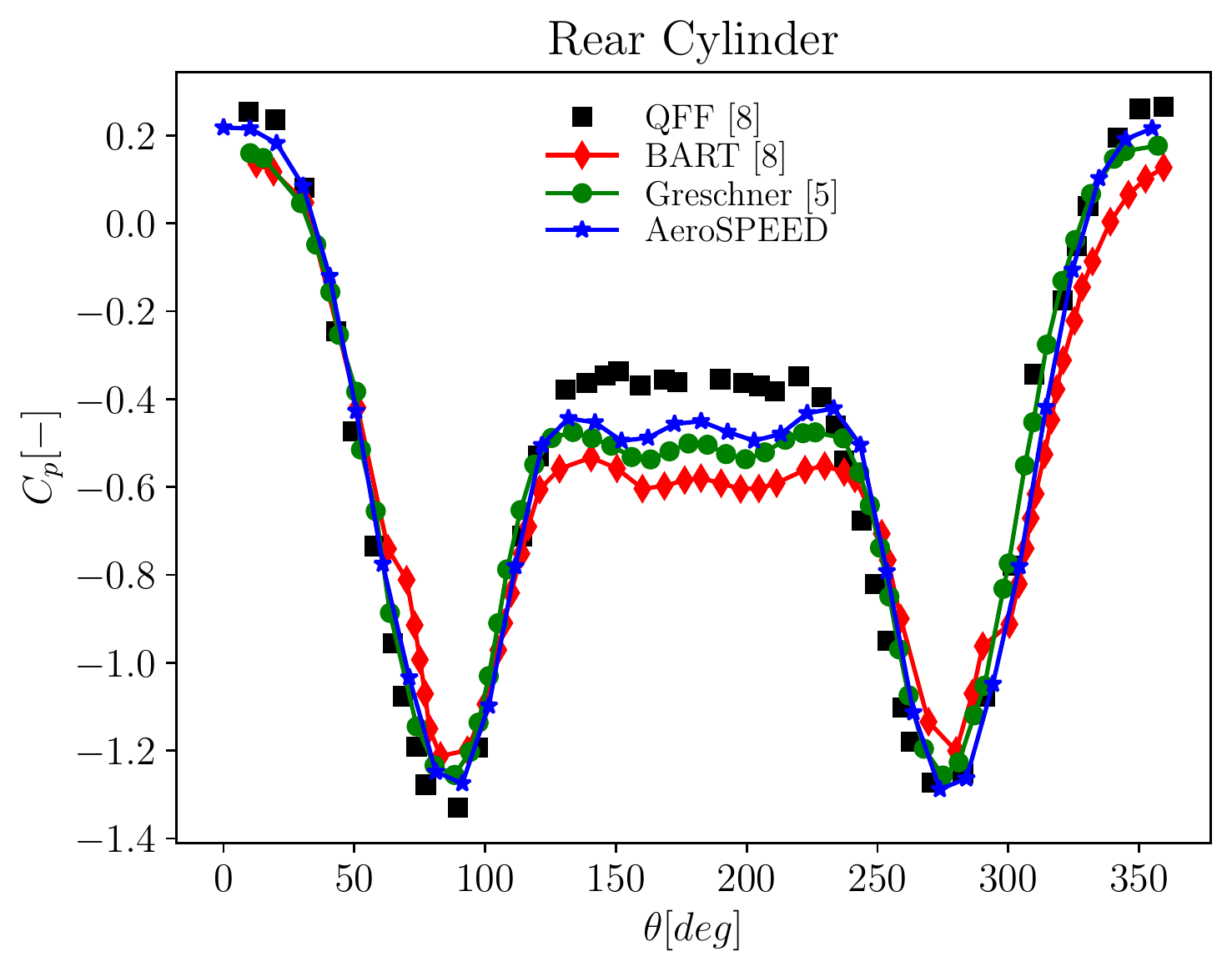}
	\end{subfigure}}
	\caption{ Average $C_P$ distribution on the surface of the front and rear cylinders, where ${\displaystyle C_P = \frac{\overline{p}-p_{ref}}{\frac{1}{2}\rho_0 U_\infty^2}}$.}
	\label{fig:CilPressure}
\end{figure}

The average flow field is characterized by a mostly symmetric recirculation regions after the cylinders, see Fig.~\ref{fig:flowfield} (left). 
The first recirculation length is about $2 D$, aligned with the results of \cite{Greschner2012}. A visualization of the vortex structures in the flow field, see Fig.~\ref{fig:flowfield} (right), is made by employing the Q criterion, where $Q = \frac{1}{2}\Big(\trace{( \nabla{\mathbf{u}})}^2 + \trace{(\nabla\mathbf{u} \nabla\mathbf{u})} \Big)$. 
In Fig.~\ref{fig:IntegralValues}, we compare the prediction on the force coefficients with the results collected in \cite{Lockard2011}. The results are quite heterogeneous due to the complexity of the problem and the numerous different strategies among the research groups. We denote with $C_D$ and $C_L$ the mean drag and lift coefficients, and their root mean squared $rms$ values as $\displaystyle \widetilde{C_{L}}, \widetilde{C_{D}}$, respectively.
All the computed integral results are within the standard deviation from the literature data. For further comparison, we consider Fig.~\ref{fig:CilPressure}, where we plot the values of the pressure coefficient $C_p$ along the cylinders and we compared it with the experiments \cite{Lockard2011}, \cite{BART} and with the computations performed in \cite{Greschner2012}.
The main differences are located in the aft of both cylinders, and are due to the prediction of a different pressure recovery, resulting in small shifts in the separation point locations.
Overall, our flow predictions are aligned with the references.

\subsubsection{Aeroacoustic Validation.}
\begin{figure}
 \adjustbox{valign=m}{
 \includegraphics[width=0.065\textwidth] {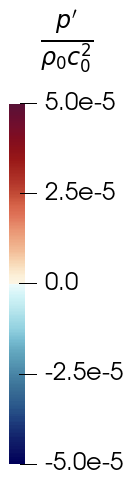}
 }
 \adjustbox{valign=m}{
	\includegraphics[width=0.39\textwidth]{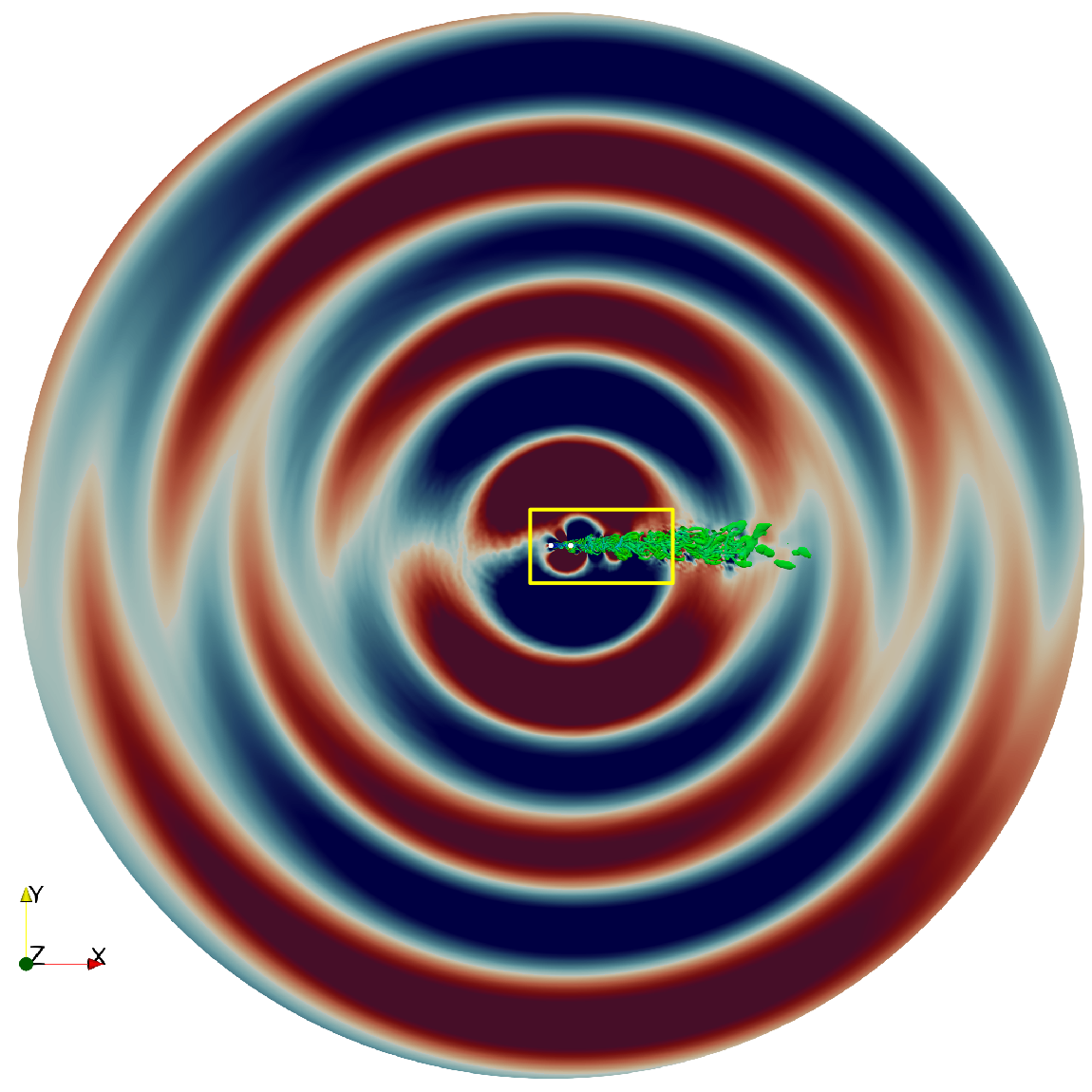}
 }
  \adjustbox{valign=m}{
 \includegraphics[width=0.39\textwidth]{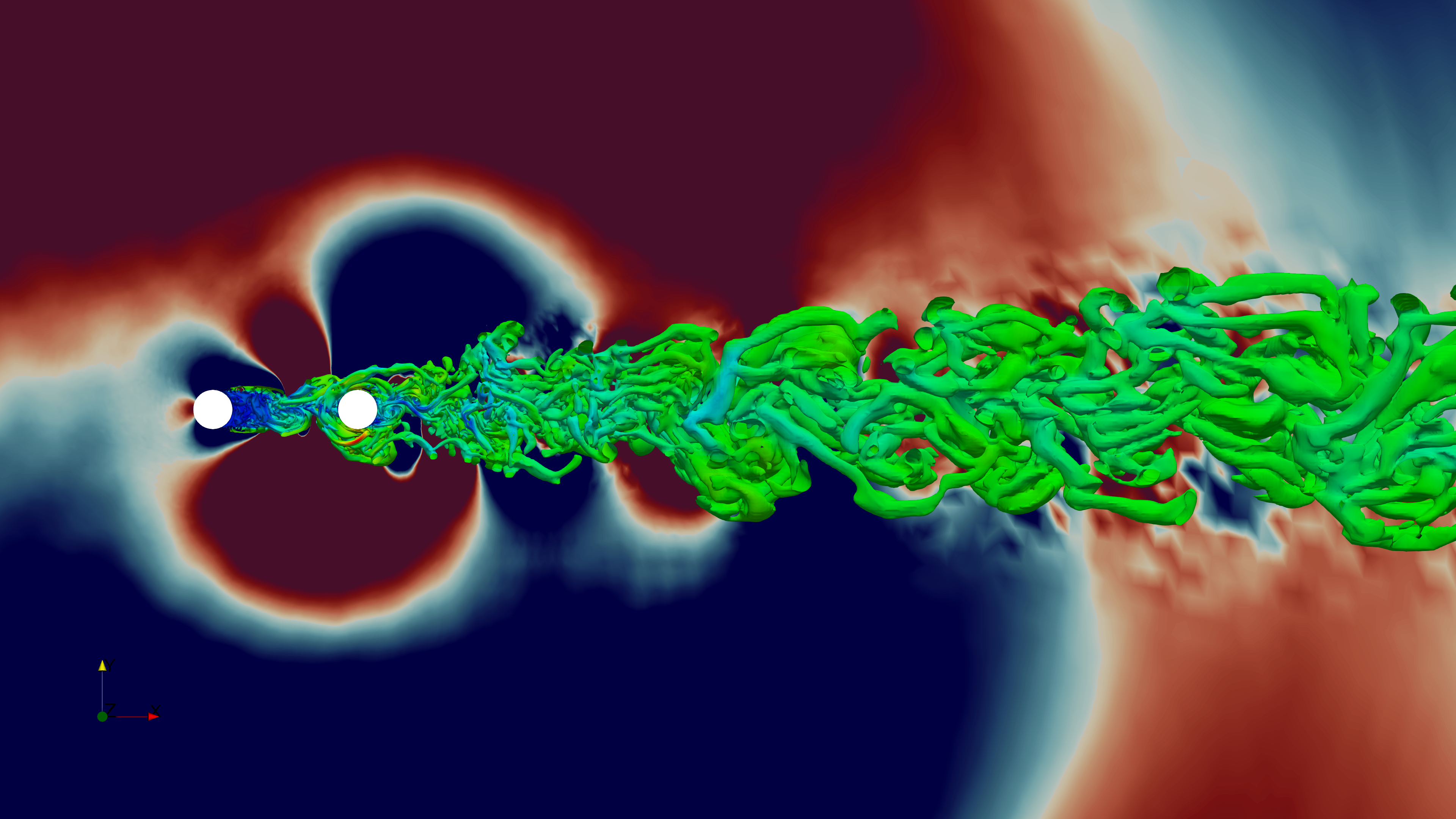}  
 }
  \adjustbox{valign=m}{
\includegraphics[width=0.065\textwidth] {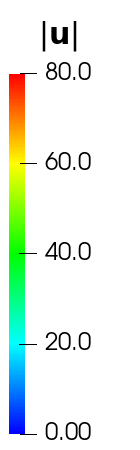} 
}
	\caption{Snapshot of the fluctuating pressure (from the acoustic computations), and Q criterion  colored with the velocity magnitude (from the flow computations).	\label{fig:SnapshotStructures}
}
\end{figure}

\begin{figure}
	\centering
	\includegraphics[valign=m,width=0.40\textwidth]{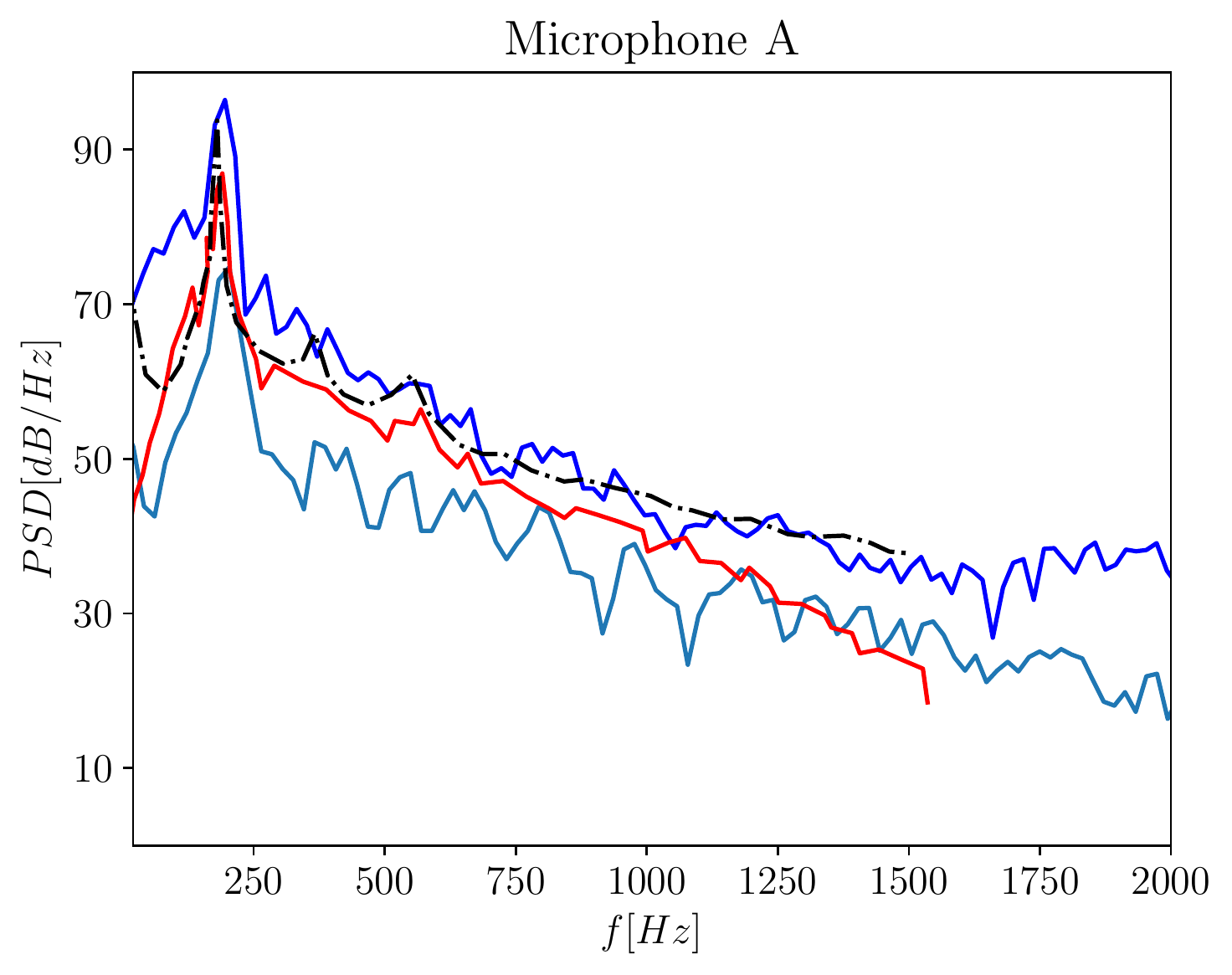} \hfill
	\includegraphics[valign=m,width=0.40
	\textwidth]{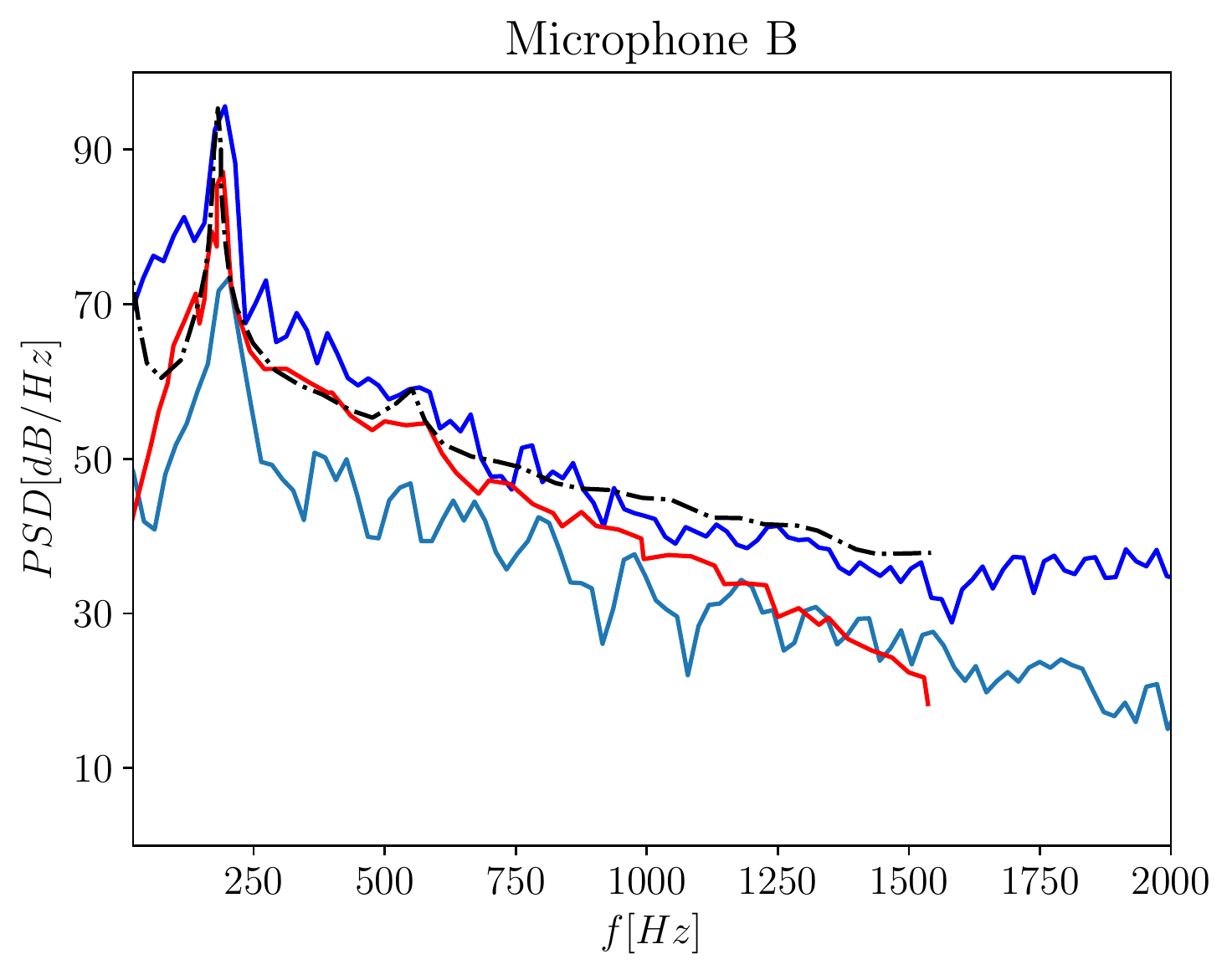}
    \includegraphics[valign=m,width=0.15
	\textwidth]{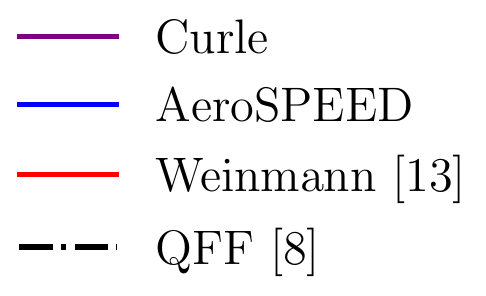}
	\caption{Comparison of the sound spectra at microphone A=(-8.33D,27.82D) and B=(9.11D, 32.49D).}
	\label{fig:SoundSPECTRA}
\end{figure}

In Fig.~\ref{fig:SnapshotStructures}, a snapshot of the pressure fluctuations induced by the tandem cylinder is shown. The acoustic pressure fluctuations are dominated by a dipole pattern induced by the lift force acting on the rear cylinder. As suggested by Fig.~\ref{fig:IntegralValues}, the main contribution to the sound generation is from the rear cylinder, being its $\widetilde{C}_{L}$ much larger than the front cylinder one.
Also, we can compare the different structures coming from the flow with the bigger structures solved by the acoustics. In Fig.~\ref{fig:SoundSPECTRA} we compare the sound spectra obtained with different methodologies, such as the Curle method described in Section \ref{sec:Curle}, another Curle analogy with a spanwise corrections proposed in \cite{Weinmann2010}, experimental data from QFF \cite{Lockard2011} and the results computed with AeroSPEED.
We observe that, although all the aeroacoustic solutions predict the peak coming from the lift frequency of the rear cylinder, AeroSPEED better matches the $PSD$ values of the experimental data.

\section{Conclusion}
We introduce AeroSPEED, a solver for aeroacoustic problems that couples a finite volume solution onto a spectral element space and solves the Lighthill's wave equation.
The high-order spectral acoustic solver was compared with the commercial software COMSOL on a model problem with a manufactured solution. The accuracy of the numerical results agree with theoretical estimates and the performance of the two solvers is compared in terms of accuracy versus computational time. 
The spectral element solution obtained with AeroSPEED is able to guarantee higher accuracy with lower computational time. 
We then applied our solver to simulate the acoustic propagation inside a simplified car cockpit (the Noise Box).
The solution obtained with the two different solvers for the pressure signals in two different locations inside the domain perfectly match. 
Next, we studied a more complex application where the noise is generated by the highly unsteady flow around tandem cylinders. We compared the results obtained with AeroSPEED with experimental and numerical tests performed by many research groups on the tandem cylinder benchmark case, proving the prediction capabilities of the proposed approach also for relevant and more challenging aeroacoustic problems.
%
%

\end{document}